\documentclass[a4paper, 12pt]{article}

\usepackage{ucs}                
\usepackage[utf8x]{inputenc} 	

\usepackage[T1]{fontenc}     
\usepackage{lmodern}         
\usepackage{xspace}

\usepackage[text={165mm,250mm},centering]{geometry}

\usepackage{rcs}  \RCS $Revision: 2.9 $ \RCS $Date: 2022/05/23 21:09:29 $ \RCS $Author: hgm $ \RCS $RCSfile: 19_low-rank_tensor_prob_den_char.tex,v $ \RCS $Id: 19_low-rank_tensor_prob_den_char.tex,v 2.9 2022/05/23 21:09:29 hgm Exp $

\usepackage{upgreek}
\usepackage{color}
\usepackage{amsmath}
\usepackage{amssymb}
\usepackage{amsbsy}
\usepackage{amsthm}
\usepackage{wasysym}
\usepackage{tocbibind}          
\usepackage{float}
\usepackage{listings,paralist}
\usepackage[ruled]{algorithm}
\usepackage{tensor}
\usepackage{graphicx}
\usepackage{caption}
\usepackage{algorithm}
\usepackage{algpseudocode}
\usepackage[dvipsnames]{xcolor}
\usepackage{authblk}

\usepackage{url}
\DeclareUrlCommand\email{\urlstyle{tt}}
\usepackage[breaklinks,colorlinks=true,plainpages=false,linkcolor=blue]{hyperref}

\usepackage[british]{babel}
\newcommand{\citep}[1]{\cite{#1}}

\usepackage{refdef}
\usepackage{ifontdef}

\usepackage{amsmath}

\newcommand{\ignore}[1]{}

\newcommand{\Lp}{\mrm{L}}

\newcommand{\Ck}{\mrm{C}}

\def\frks#1#2{{\raise0.8ex\hbox{{\leavevmode$\scriptstyle #1$}}}{\kern-0.4ex\hbox{/}}{\raise-0.2ex\hbox{\kern-0.4ex\hbox{{\leavevmode$\scriptstyle #2$}}}}}

\newcommand{\vrho}{\varrho}

\newcommand{\vphi}{\varphi}

\newcommand{\vek}[1]{\mathchoice{\displaystyle\boldsymbol{#1}}
{\textstyle\boldsymbol{#1}}{\scriptstyle\boldsymbol{#1}}
{\scriptscriptstyle\boldsymbol{#1}}}

\newcommand{\ops}[1]{\mathchoice{\displaystyle\mathsf{#1}}
{\textstyle\mathsf{#1}}{\scriptstyle\mathsf{#1}}
{\scriptscriptstyle\mathsf{#1}}}

\newcommand{\tnb}[1]{\mathchoice{\displaystyle\mathboldsans{#1}}
{\textstyle\mathboldsans{#1}}{\scriptstyle\mathboldsans{#1}}
{\scriptscriptstyle\mathboldsans{#1}}}

\newcommand{\tns}[1]{\mathchoice{\displaystyle\mathsans{#1}}
{\textstyle\mathsans{#1}}{\scriptstyle\mathsans{#1}}
{\scriptscriptstyle\mathsans{#1}}}

\newcommand{\vbar}[1]{\vek{\bar{#1}}}
\newcommand{\vtil}[1]{\vek{\tilde{#1}}}
\newcommand{\vhat}[1]{\vek{\hat{#1}}}

\newcommand{\Ttil}[1]{\tnb{\tilde{#1}}}
\newcommand{\That}[1]{\tnb{\hat{#1}}}

\newcommand{\EXP}[1]{\mathbb{E}\left(#1\right)}

\DeclareMathOperator{\tr}{tr}

\DeclareMathOperator{\supp}{supp}
\DeclareMathOperator{\sign}{sign}
\DeclareMathOperator{\cl}{cl}

\newcommand{\dd}{\mathop{}\!\partial}
\newcommand{\di}{\mathop{}\!\mathrm{d}}
\newcommand{\Di}{\mathop{}\!\mathrm{D}}
\newcommand{\ii}{\mathchoice{\displaystyle\mathrm i}
{\textstyle\mathrm i}{\scriptstyle\mathrm i}
{\scriptscriptstyle\mathrm i}}
\newcommand{\ee}{\mathchoice{\displaystyle\mathrm e}
{\textstyle\mathrm e}{\scriptstyle\mathrm e}
{\scriptscriptstyle\mathrm e}}

\newcommand{\KL}{Karhunen-Lo\`eve}
\newcommand{\ip}[2]{\left\langle #1 , #2 \right\rangle}
\newcommand{\bkt}[2]{\left\langle #1 \mid #2 \right\rangle}

\newcommand{\ns}[1]{\left| #1 \right|}
\newcommand{\nd}[1]{\left\Vert #1 \right\Vert}

\newcommand{\trpos}{{\ops{T}}}

\newcommand{\bbbone}{{\mathchoice {\rm 1\mskip-4mu l} {\rm 1\mskip-4mu l}
{\rm 1\mskip-4.5mu l} {\rm 1\mskip-5mu l}}}

\newcommand{\R}{\D{R}}
\newcommand{\N}{\D{N}}
\newcommand{\bbC}{\D{C}}
\newcommand{\Tp}{\C{T}}
\DeclareMathOperator{\cov}{cov}

\newcommand    {\RR}
                {\D{R}}

\newcommand    {\ten}    { \otimes }
\newcommand    {\Ten}    { \bigotimes }
   \newcommand    {\Rd}    { \D{R}^d }
   \newcommand    {\Rn}    { \D{R}^n }
   \newcommand    {\RN}    { \D{R}^N }
    
    \newcommand    {\Rm}    { \D{R}^m }

\newcommand{\alphab}{\vek{\alpha}}
\newcommand{\betab}{\vek{\beta}}

\newcommand{\frkt}[2]{{{\raise0.6ex\hbox{{\leavevmode$\textstyle #1$}}}{\raise0.25ex\hbox{\kern-0.35ex\hbox{/}}}{\raise-0.3ex\hbox{\kern-0.4ex\hbox{{\leavevmode$\textstyle  #2$}}}}}}
\newcommand{\frkc}[2]{{{\raise0.7ex\hbox{{\leavevmode$\scriptstyle #1$}}}{\raise0.2ex\hbox{\kern-0.4ex\hbox{\footnotesize /}}}{\raise-0.2ex\hbox{\kern-0.4ex\hbox{{\leavevmode$\scriptstyle #2$}}}}}}
\newcommand{\frkx}[2]{{{\raise0.75ex\hbox{{\leavevmode$\scriptscriptstyle #1$}}}{\raise0.17ex\hbox{\kern-0.45ex\hbox{\scriptsize /}}}{\raise-0.15ex\hbox{\kern-0.4ex\hbox{{\leavevmode$\scriptscriptstyle #2$}}}}}}
\newcommand{\frkz}[2]{{{\raise0.75ex\hbox{{\leavevmode$\scriptscriptstyle #1$}}}{\raise0.17ex\hbox{\kern-0.45ex\hbox{\tiny /}}}{\raise-0.15ex\hbox{\kern-0.4ex\hbox{{\leavevmode$\scriptscriptstyle #2$}}}}}}

\newcommand{\frk}[2]{{\mathchoice{{\frkt{#1}{#2}}}{{\frkc{#1}{#2}}}{{\frkx{#1}{#2}}}{{\frkz{#1}{#2}}}}}

\newcommand{\bSigma}{\vek{\Sigma}}

\newcommand{\Ip}{\C{I}}
\newcommand{\nin}{\notin}

\newcommand{\mybold}[1]{\vek{#1}}

\newcommand{\bx}{\mybold{x}}
\newcommand{\by}{\mybold{y}}

\newcommand{\bt}{\mybold{t}}

\newcommand{\bC}{\mybold{C}}
\newcommand{\bI}{\mybold{I}}

\newcommand{\J}{\C{J}}

\newcommand{\vX}{\vek{\xi}}

\newcommand{\thetab}{\vek{\theta}}
\newcommand{\Fd}{\C{F}_d}
\newcommand{\iFd}{\C{F}_{d}^{-1}}

\newcommand{\etab}{\vek{\eta}}

\newcommand{\cdf}{\texttt{cdf}\xspace} 
\newcommand{\pdf}{\texttt{pdf}\xspace} 
\newcommand{\pdfs}{\texttt{pdfs}\xspace} 
\newcommand{\pcfs}{\texttt{pcfs}\xspace} 
\newcommand{\pcf}{\texttt{pcf}\xspace} 

\newcommand{\BIGOP}[1]{\mathop{\mathchoice%
{\raise-0.22em\hbox{\huge $#1$}} {\raise-0.05em\hbox{\Large $#1$}}
{\hbox{\large $#1$}}{#1}}}
\newcommand{\BIGboxplus}{\mathop{\mathchoice%
{\raise-0.35em\hbox{\huge $\boxplus$}}%
{\raise-0.15em\hbox{\Large $\boxplus$}}{\hbox{\large
$\boxplus$}}{\boxplus}}}
\newcommand{\bigtimes}{\BIGOP{\times}}

\newtheorem{theorem}{Theorem}[section]
\newtheorem{remark}[theorem]{Remark}
\newtheorem{definition}[theorem]{Definition}
\newtheorem{example}[theorem]{Example}

\newcommand{\bz}{\mybold{z}}
\newcommand{\bs}{\mybold{s}}

\definecolor{light-blue}{rgb}{0.1,0.75,0.9} 
\definecolor{myred}{rgb}{0.8, 0.1, 0.2}

\newcommand{\textdate}{}
\date{\textdate}

\newcommand{\tand}{, }
\newcommand{\nand}{$\cdot\;$}
\newcommand{\autheadcr}{\\}
\newcommand{\thetitle}{Computing f-Divergences and Distances of High-Dimensional Probability Density Functions}
\newcommand{\thesubtitle}{Low-Rank Tensor Approximations}

\newcommand{\authorAL}{Alexander~Litvinenko}

\newcommand{\authorYM}{Youssef~Marzouk}

\newcommand{\authorHGM}{Hermann~G.~Matthies}

\newcommand{\authorMS}{Marco~Scavino}

\newcommand{\authorAS}{Alessio~Spantini}

\newcommand{\affilTUBS}{\small 
                        Technische Universit\"at Braunschweig, Brunswick, Germany}
\newcommand{\affilMIT}{\small MIT, Cambridge (MA), USA}
\newcommand{\affilKAUST}{\small Universidad de la Rep\'ublica, Instituto de Estad\'{\i}stica (IESTA), Montevideo, Uruguay}
\newcommand{\affilRWTH}{\small Rheinisch-Westf\"alische Technische Hochschule (RWTH)
                         Aachen, Germany}

\newcommand{\thesubject}{41A05 \nand 41A45 \nand 41A63 \nand 60-08 \nand 60E07 \nand 60E10 \nand 62-08 
                         \nand 62E17 \nand 62H10 \nand 65C20 \nand 65F55 \nand 65F60 \nand 94A17}

\newcommand{\thekeywords}{high-dimensional probability density\tand Kullback-Leibler divergence\tand
        f-divergence \tand distance \tand tensor representation\tand computational algorithms \tand
         low-rank approximation}

\begin{document}

\title{\thetitle{}
    \autheadcr      --- \thesubtitle{} ---
 }


%

\urldef{\ALemail}{\email}{Litvinenko@uq.rwth-aachen.de}
\author[a]{\authorAL
\thanks{Corresponding author: RWTH Aachen, 52072 Aachen, Germany,  \autheadcr 
{e-mail: }{\ALemail}
}
}
\author[b]{\authorYM}
\author[c]{\authorHGM}
\author[d]{\authorMS}
\author[b]{\authorAS}

\affil[a]{\affilRWTH}
\affil[b]{\affilMIT}
\affil[c]{\affilTUBS}
\affil[d]{\affilKAUST}

%

\maketitle




\begin{abstract}
Very often, in the course of uncertainty quantification tasks or 
data analysis, one has to deal with high-dimensional random variables.
Here the interest is mainly to compute characterisations like the entropy, 
the Kullback-Leibler divergence, more general $f$-divergences, or other such characteristics based on
the probability density.  The density is often not available directly, 
and it is a computational challenge to just represent it in a numerically
feasible fashion in case the dimension  is even moderately large.   It
is an even stronger numerical challenge to then actually compute said characteristics
in the high-dimensional case.  
In this regard
 it is proposed
to approximate the discretised density in a compressed form, in 
particular by a low-rank tensor.
This can alternatively be obtained from the corresponding 
probability characteristic function,
or more general representations of the underlying random variable.
The mentioned characterisations need 
point-wise functions like the logarithm.  
This normally rather trivial task becomes computationally difficult
when the density is approximated in a compressed resp.\ low-rank tensor 
format, as the point values are not directly accessible.
The computations become possible by considering the compressed data 
as an element of an associative, commutative algebra with an inner product,
and using matrix algorithms to accomplish the mentioned tasks.
The representation as a low-rank element of a high
order tensor space allows to reduce the computational complexity and storage cost from 
exponential in the dimension to almost linear.  


\vspace{1em}
{\noindent\textbf{Keywords:} \thekeywords}

\vspace{1em}
{\noindent\textbf{MSC Classification:} \thesubject}

\end{abstract}

%
%
%
%
%
%
%
%
%
%
%
%
%
%












\section{Introduction}  \label{S:intro}
%


In statistics and probability, and in particular in the vigorous field of
uncertainty quantification (UQ) \citep{HandbookUQ17, Matthies_encicl, matthiesKeese05cmame}, 
one often has to deal with high-dimensional random variables (RVs) with values in $\Rd$. 
%
Frequently, one is interested in characteristic quantities of interest (QoI), examples are
the differential entropy, the Kullback-Leibler  or more general $f$-divergences, 
or other such characteristics based on integrals of some function of 
the probability density function (\pdf).  
Computing such integrals, or just storing a discretised form of the \pdf, may become a
daunting task for even moderately high dimensions $d$.

\subsection{Motivation and main idea} \label{SS:motiv}
One apparently quite efficient way to deal with such high-dimensional
RVs is to view them as elements of some tensor product space of tensors of
high order --- typically the order corresponds to the dimension $d$ \citep{hackbusch2012tensor}.  
Such tensors can be approximated by low-rank tensors, and thus it becomes feasible to
deal with them numerically. 
Here we concentrate solely
on real valued continuous RVs which possess a \pdf,
so another possibility is the
representation of the \pdf in some computationally advantageous form.

Possible QoIs are often
defined or only accessible in some specific representation of the RV.
We are thinking of QoIs such as the
entropy, the Kullback-Leibler (KL) divergence, or more generally $f$-divergences, 
as well as
Hellinger distances, and (central) and generalised moments, to name a few such QoIs
\citep{LieseVajda06, NielsenNock2013, Toulias2020}. 
Some are defined and can efficiently be computed with
a functional representation of the RV, and some are defined and efficiently computed
if one has access to the point values of the \pdf, the case to be considered here.

The algorithms we formulate are for a compressed discrete representation of the \pdf, and
we actually use for our computations a low-rank tensor approximation.
We sketch three possibilities how to arrive at such a low-rank
representation of the density.  One is via tensor function representations 
\citep{Beylkin05, bebe-aca-2011, hackbusch2012tensor, 
GoroMarzouk2019, BigoniEtal2016} of the \pdf \citep{DolgovFox19, Dolgov20},
another via an analogous representation of the probabilistic characteristic function
(\pcf) \citep{Shephard, Zoltan13, Witkovsk2016NumericalIO}, and a third one
is via a low-rank representation of the RV itself---e.g.\ 
\citep{dolgov2014computation, dolgov2015polynomial}.

Connecting and connected with these different representations---next to the 
already mentioned \pcf ---are other well known characterising objects like the 
moment generating function or
the second characteristic or cumulant generating function, and we show how these may be
efficiently computed as well in a low-rank tensor format.
Motivating factors for an accuracy controlled low-rank tensor compression include the following
 \citep{Kolda:07, hackbusch2012tensor}:
\begin{compactitem}
\item reduce the storage cost to linear or even sub-linear in the dimension $d$;
\item perform algebraic operations in cost linear or even sub-linear in the dimension $d$;
\item combining it with the fast Fourier transform (FFT) yields a \emph{superfast} 
   FT \citep{nowak2013kriging, DoKhSa-qtt_fft-2011}.
\end{compactitem}
On the other hand, general limitations of such a tensor compression technique are that 
\begin{compactitem}
\item it could be time-consuming to compute 
     a low-rank tensor decomposition; 
\item with sampling or point evaluation, it requires an axes-parallel mesh; 
\item after algebraic manipulations, a re-compression may be necessary;
\end{compactitem}
It is a fact that there are situations where storage of
all items is not feasible, or operating on all items is out of the question,
so that some kind of compression has to be employed, and
low-rank tensor techniques offer a very promising avenue \citep{Kolda:07, hackbusch2012tensor}.
But for the algorithms shown later in \refS{algs} all that is required is that the compressed
data can be considered
as an element of an associative, commutative algebra with an inner product.

Further, the actual computational representation is allowed to be a lossy compression, and 
the algebraic operations may be performed in an approximate fashion, so
as to maintain a high compression level.  We address explicitly
the representation of data as a tensor with compression in the form of a 
low-rank representation.
The suggested technique, in order to be employed efficiently, assumes that all involved 
tensors have a low rank, and the rank is not increasing strongly after linear algebra operations.
In \refS{Numerics} we show some exemplary numerical examples of 
applications of this framework to approximate high-dimensional \pdfs of
$\alpha$-stable type---which are given only through their \pcf---and compute 
divergences and distances between such probability distributions.

Let us give an example which motivates much of the following formulation and development.
Assume that one is dealing with a random vector 
$\vX$ in a high-dimensional vector space $\Rd$, i.e.\
$   \vX=( \xi_1,...,\xi_d )^\trpos : \Omega\to\Rd $,
defined on a suitable probability space $\Omega$.  Further assume that this random vector
has a \pdf $p_{\vX}: \Rd\to\D{R}$, and that we
want to compute---as a simple example of the possible tasks envisioned---the
\emph{differential entropy} (see also \refeq{eq:discr-diff-entrop}), the
\emph{expectation} of the negative logarithm of the \pdf.  This requires the 
point-wise logarithm of $p_{\vX}$:
\begin{equation}  \label{eq:exmpl-diff-entrop}
   h(p_{\vX}) := \EXP{-\ln(p_{\vX}(\by))} := \int_{\Rd} -\ln(p_{\vX}(\by)) p_{\vX}(\by)\, \di\by .
\end{equation}
Even if one has an analytical expression resp.\ approximation
for the \pdf $p_{\vX}$, it still may not be
possible to compute the above integral analytically, and so we propose to do this
numerically.  This immediately provokes the \emph{curse of dimensionality}, as the
numerical evaluation of such high dimensional integrals can be very expensive, with
work proportional $n^d$, where $n\in\D{N}$ is a discretisation parameter .
In the following, we shall propose one
possible approach to alleviate the computational burden.

The goal is a discrete low-rank tensor point evaluation of the \pdf.  The first
starting point 
is the assumption that one has or may obtain
a low-rank tensor function representation 
\citep{Beylkin05, bebe-aca-2011, hackbusch2012tensor,
GoroMarzouk2019, BigoniEtal2016} of the \pdf
\citep{DolgovFox19, Dolgov20}.  The second possible starting point in 
is that one has such a low-rank tensor function of the \pcf
\citep{Lukacs70b, Shephard, Zoltan13, Witkovsk2016NumericalIO, belomIosipoi2021}.
Then via the Fourier transform one may come back to the \pdf.  And the third possibility
sketched 
is that one has a sparse or low-rank representation of a 
high-dimensional RV.  From this one may evaluate the \pcf on a
tensor grid, and from there then via the Fourier transform again arrive at a low-rank
approximation of the \pdf.
%
%
%

%
Further assume that the \pdf 
of the high-dimensional RV
has its support in a compact hyper-rectangle,
 $ \supp\, p_{\vX} := \cl \{ \by\in\Rd \mid p_{\vX}(\by) \ne \vek{0} \} \subseteq
\bigtimes_{\nu=1}^d [\xi_\nu^{(\min)}, \xi_\nu^{(\max)}] \subset \Rd $.
A fully discrete representation of the \pdf can be achieved (further details will be
provided in \refS{theory}) by picking
in each dimension $1\le \nu \le d$ of $\Rd$ an equidistant grid vector $\vhat{x}_\nu :=
(\hat{x}_{1,\nu},\dots,\hat{x}_{M_{\nu},\nu})$ of size $M_\nu$, such that for all $\nu$ it holds that
$\hat{x}_{i_\nu,\nu} \in [\xi_\nu^{(\min)}, \xi_\nu^{(\max)}]$ for $1\le i_\nu \le M_\nu$ .
The size $M_{\nu}$ could be different for each
dimension $\nu$, but for the sake of simplicity here we shall often assume them all equal to $n$,
i.e.\ each $\vhat{x}_{\nu} \in \Rn$.  The whole grid will be denoted by 
$\That{X} = \bigtimes_{\nu=1}^d \vhat{x}_{\nu}  = (\hat{\tns{X}})_{(\nu,i_1,\dots,i_d)}$
with $1\le i_\nu \le M_\nu$.

The notation $\tnb{P}:=p_{\vX}(\That{X})$ will denote the tensor $\tnb{P}\in 
\bigotimes_{\nu=1}^d \R^{M_{\nu}} =: \C{T} = (\Rn)^d = \RN$, with
$\dim\C{T} = \prod_{\nu=1}^d M_{\nu} =: N = n^d $, the components of which are
the evaluation of the \pdf $p_{\vX}$ on the grid $\That{X}$
\begin{equation}   \label{eq:pdf-tensor-1}
   \tnb{P}:=p_{\vX}(\That{X}) := (\tns{P}_{i_1,\dots,i_d}) := 
   (p_{\vX}(\hat{x}_{i_1,1},\dots,\hat{x}_{i_d,d})) ,
   \quad 1 \le i_\nu \le M_\nu,\; 1\le\nu\le d .
\end{equation}
If $n$ and especially the dimension $d$ are even moderately large, the total dimension
$N=n^d$ is very quickly a huge number, maybe even so that it is not possible to
reasonably store that amount of information.

Sometimes the density $p_{\vX}$ is---at least approximately---given as an analytical
expression, and it may be possible to approximate it \citep{HA12} through
a low-rank function tensor representation.  In the simplest case this would look like
\begin{equation}
\label{eq:p}
p_{\vX}(\by) \approx \tilde{p}_{\vX}(\by)=\sum_{\ell=1}^R \Ten_{\nu=1}^d p_{\ell,\nu}(y_\nu),
\end{equation}
where each $p_{\ell,\nu}$ is only a function of the real variable $y_\nu$ in dimension $\nu$.
This so-called \emph{canonical polyadic} (CP) tensor representation \citep{HA12}
(cf.\ \refS{tensor-rep}) becomes
computationally viable when the \emph{rank} $R$ can be chosen fairly small.   
The multi-dimensional Gaussian distribution with diagonal covariance matrix is an
obvious simple case in point with $R=1$. 
The $p_{\ell,\nu}$ are evaluated on the grid vector $\vhat{x}_{\nu}$
for all $\nu$ and $\ell$,  giving
$
\vek{p}_{\ell,\nu} := (p_{\ell,\nu}(\hat{x}_{1,\nu}),\dots,p_{\ell,\nu}(\hat{x}_{M_{\nu},\nu})) 
\in \D{R}^{M_\nu} .
$
This is a building block for a possible low-rank CP representation of the tensor $\tnb{P}$,
as now
\begin{equation}  \label{eq:bsc-CP-repr}
 \tnb{P}\approx \sum_{\ell=1}^R \Ten_{\nu=1}^d \vek{p}_{\ell,\nu}.
\end{equation}
More details and descriptions about this and other low-rank tensor approximations will be
given in \refS{tensor-rep}.
To evaluate now numerically an expression like the differential entropy
\refeq{eq:exmpl-diff-entrop}, the integral is replaced by a numerical quadrature
\begin{equation}  \label{eq:exmpl-diff-entrop-discr}
  h(p_{\vX}) \approx \sum_{i_1=1}^{M_1} \cdots \sum_{i_d=1}^{M_d}  -\ln(\tns{P}_{i_1,\dots,i_d})
   \tns{P}_{i_1,\dots,i_d} w_{i_1,\dots,i_d},
\end{equation}
where $w_{i_1,\dots,i_d}$ are integration weights, which will all be chosen equal
$w_{i_1,\dots,i_d} \propto N^{-1}$, cf.\ \refS{theory} and \refS{statistics}.

The challenge in \refeq{eq:exmpl-diff-entrop-discr}
is the huge number of terms in the sum.  Here the low-rank representation \refeq{eq:bsc-CP-repr}
can be used to advantage, cf.\ \refS{theory}, but then the challenge is to compute
the logarithm in \refeq{eq:exmpl-diff-entrop-discr} when $\tnb{P}$ is in 
the representation \refeq{eq:bsc-CP-repr}. 
For this certain algebraic properties of the
space $\C{T}$ will be used, as will be explained in \refS{theory} and \refS{algs}.
There are of course other ways to arrive at a discrete low-rank representation of a
high-dimensional function.  One group of such possible methods are the various 
``cross''-procedures \citep{EspigGrasedyckHackbusch2009, oseledetsTyrt2010, 
bebe-aca-2011, BallaniGrasedyck2015, ds-alscross-2017pre, dolgovEtal2020}.  In any case, this is
not supposed to be an exhaustive survey of such methods, and only intended
to give some hints and point to some possibilities.

%
At other times the \pdf may not be available, but instead the probabilistic characteristic
function (\pcf) $\vphi_{\vX}$ of the \pdf 
\begin{equation}  \label{eq:def-char-fct}
  \vphi_{\vX}(\bt):=\EXP{\exp(\ii \bkt{\vX}{\bt})} :=
    \int_{\D{R}^d} p_{\vX}(\by)\exp(\ii \bkt{\by}{\bt})\,\di\by =: \Fd(p_{\vX})(\bt), 
\end{equation}
where $\bt=(t_1,t_2,...,t_d) \in \D{R}^d$ is the dual variable to $\by\in\Rd$,
$\bkt{\by}{\bt}=\sum_{j=1}^d y_jt_j$ is the canonical inner product on $\Rd$, and
$\C{F}_d$ is the probabilist's $d$-dimensional Fourier transform (FT)  \refeq{eq:def-char-fct},
may be given \citep{belomIosipoi2021}, at least 
approximately, as an analytical expression
(similar to \refeq{eq:p}) 
\begin{equation}  \label{eq:pcf2}
 \vphi_{\vX}(\bt)\approx  \tilde{\vphi}_{\vX}(\bt)=\sum_{\ell=1}^R \Ten_{\nu=1}^d 
 \vphi_{\ell,\nu}(t_\nu),
\end{equation}
where the $\vphi_{\ell,\nu}(t_\nu)$ are one-dimensional functions.  Examples of
such a situation are \emph{elliptically contoured $\alpha$-stable distributions}, or also
\emph{symmetric infinitely divisible distributions} \citep{belomIosipoi2021}.
From this it may be deduced that an approximate low-rank expression of the \pdf
is given by  (cf.\ \refeq{eq:p})
\begin{equation}  \label{eq:motiv_pdf_lr}
  p_{\vX}(\by) \approx \tilde{p}_{\vX}(\by)= \iFd(\vphi_{\vX}){\by}
  =\sum_{\ell=1}^R \Ten_{\nu=1}^d \C{F}^{-1}_1(\vphi_{\ell,\nu})(y_\nu),
\end{equation}
where $\C{F}^{-1}_1$ is the \emph{one-dimensional} probabilist's inverse FT.

For the discrete grid $\That{X}$, there is a corresponding \emph{dual grid} 
$\That{T}= (\hat{\tns{T}})_{(\nu,i_1,\dots,i_d)}$ with $1\le i_\nu \le M_\nu$
for the discrete Fourier transform \citep{bracewell} 
of same size and dimensions.  Similarly to the \pdf evaluated on the grid $\That{X}$
and represented by the tensor $\tnb{P}$, 
the \pcf will be used in a discrete setting evaluated on this dual grid $\That{T}$
\begin{equation}  \label{eq:grid-char}
   \tnb{\Phi}:=\vphi_{\vX}(\That{T}) := 
   (\vphi_{\vX}(\hat{t}_{i_1,1},\dots,\hat{t}_{i_d,d}))=\Fd(\tnb{P}) ,
   \quad 1 \le i_\nu \le M_\nu,\; 1\le\nu\le d .
\end{equation}
From \refeq{eq:pcf2} it is now easy to see that, as in \refeq{eq:bsc-CP-repr},
$
  \tnb{\Phi} \approx \sum_{\ell=1}^R \Ten_{\nu=1}^d \vek{\vphi}_{\ell,\nu},
$
where the vectors $\vek{\vphi}_{\ell,\nu}$ are the evaluations of the one-dimensional
functions $\vphi_{\ell,\nu}$ from \refeq{eq:pcf2}, given by
$
   \vek{\vphi}_{\ell,\nu} := (\vphi_{\ell,\nu}(\hat{t}_{1,\nu}),\dots, 
     \vphi_{\ell,\nu}(\hat{t}_{M_{\nu},\nu})) \in \D{R}^{M_\nu} .
$
From these vectors $\vek{\vphi}_{\ell,\nu}$ one may now compute with the discrete one-dimensional
inverse Fourier transform, for the sake of simplicity again denoted by $\C{F}^{-1}_1$,
the vectors $\vek{p}_{\ell,\nu} = \C{F}^{-1}_1(\vek{\vphi}_{\ell,\nu})$, such that from 
\refeq{eq:motiv_pdf_lr} one arrives at an expression corresponding  to \refeq{eq:bsc-CP-repr}:
\begin{equation}  \label{eq:grid-char-FT}
 \tnb{P} = \iFd(\tnb{\Phi}) \approx \sum_{\ell=1}^R \Ten_{\nu=1}^d 
 \C{F}^{-1}_1(\vek{\vphi}_{\ell,\nu})  = \sum_{\ell=1}^R \Ten_{\nu=1}^d \vek{p}_{\ell,\nu} .
\end{equation}
Thus one again obtains a numerical low-rank representation for the density
as in \refeq{eq:bsc-CP-repr}.  Obviously, the cross-methods 
\citep{EspigGrasedyckHackbusch2009, oseledetsTyrt2010, bebe-aca-2011,
BallaniGrasedyck2015, ds-alscross-2017pre, dolgovEtal2020} alluded
to there can be used here too, to directly obtain a
low-rank representation of the \pcf.  Then again, with the help of the discrete
inverse Fourier transforms, one arrives at a low-rank representation of the \pdf.


%
Random vectors of high dimension $\vX$ occur also when random fields resp.\ stochastic
processes are discretised, often given through their
so-called \KL{} expansion (KLE)
\citep{Loeve,Karhunen,Karhunen_1947_KL_Expansion}  of $\vX$, truncated at $r\in\D{N}$:
\begin{equation}
\label{eq:KLE-xi0}
  {\vX}(\omega) = \sum_{\ell=0}^{r} \lambda_\ell^{1/2}\, \zeta_\ell(\omega) \vek{v}_\ell.
\end{equation}
This places the random vector $\vX(\omega) = [\dots,
\xi_k(\omega),\dots]$, a function of the two variables $(\omega,k)$, in the
tensor product $\Lp_2(\Omega)\otimes\Rd$; and the \KL{} expansion \refeq{eq:KLE-xi0}
(another form of singular value decomposition)
is a \emph{separated} representation of this tensor of second order.  
Often the singular values $\lambda_\ell^{1/2}$
decay quickly as $\ell$ grows, so that one may obtain a good approximation
with only $r$ terms. 

To see in a nutshell where this leads to, assume further that the uncorrelated
RVs $\zeta_\ell(\omega) = \sum_{\alphab} \zeta_\ell^{(\alphab)} \Psi_{\alphab}(\thetab(\omega))$ 
may be expanded in Wiener's polynomial chaos expansion (PCE),
see e.g.\ \citep{janson97},
with multi-variate polynomials $\Psi_{\alphab}(\thetab(\omega)):= 
\prod_{j=1}^d \psi_{\alpha_j}(\theta_j(\omega))$ in \emph{iid}
standard normalised Gaussians $\thetab(\omega) = (\theta_1(\omega),\dots,\theta_d(\omega))$,
where $\alphab=(\alpha_1,\dots,\alpha_d)$ is a multi-index and the 
$\psi_{\alpha_j}(\theta_j(\omega))$ are uni-variate polynomials.
Inserting this into the KLE, one obtains a combined truncated KLE / PCE 
 $ \vX(\omega) \approx 
  \sum_{\alphab} \vX^{(\alphab)} \Psi_{\alphab}(\thetab(\omega)) $,
where $\vX^{(\alphab)}= \sum_{\ell=0}^r \lambda_\ell^{1/2}\,\zeta_\ell^{(\alphab)} \vek{v}_\ell\in\Rd$.
Arguably, now the tensor---a multi-dimensional array---
$\tnb{Z} = (\xi_k^{(\alpha_1,\dots,\alpha_d)})$ represents the RV $\vX$.
Such a tensor can---as before---be represented in the \emph{canonical polyadic} (CP) 
format \citep{HA12}, say with low \emph{CP-rank} $R$:
 $ \tnb{Z} \approx \sum_{r=1}^R \Ten_{k=0}^d \vek{z}^{(k)}_r = \sum_{r=1}^R \tnb{z}_r $.
%
%

%
Having $\tnb{Z}$ in a low-rank representation then leads to a 
low-rank representation of the random vector $\vX(\thetab)$,
to obtain a formula for a quick evaluation of the  \pcf of $\vX(\thetab)$:
 $ \vphi_{\vX}(\bt) = \EXP{\exp( \ii \bkt{\bt}{\vX(\thetab)}_{\D{R}^d})} $, 
yielding the \pcf of $\vX(\thetab)$ in a low-rank tensor function format.
Now we are back to the situation as before when the \pcf was assumed given.


\subsection{Literature review and outline of the paper}  \label{SS:lit-review}
{To compute quantities of interest (QoIs) such as the 
$f$-divergences
---see \refS{statistics} for the distances and divergences considered---functions such
as  the square root, the logarithm, etc., have to be computed point-wise
on the \pdf.  But in a compressed format 
the point values are not directly accessible.  
An important contribution of this paper is to show that such computations are still
possible efficiently by operating directly on the compressed format.  This will be enabled
by identifying both the \pcf and the \pdf as elements of algebras.  A similar idea
was already used in \citep{ESHALIMA11_1, Matthies2020}
for post-processing low-rank representations of RVs, and
here it is  extended to densities and characteristic functions.  It is well-known
that algebras can be represented as linear operators resp.\ matrices in the finite-dimensional
case (e.g.\ \citep{Segal1978}), and one can employ the spectral calculus of linear 
operators in order to obtain these point-wise evaluations in a low-rank format.
In the concrete case here this boils down to using algorithms which were developed to
compute functions of matrices \citep{NHigham}.}

Therefore, special compressed / low-rank data structures are needed.  
The CP format used above for illustrative purposes may be the simplest, but
we suggest to use the low-rank 
tensor train (TT) data format \citep{oseledetsTyrt2010, oseledets2011}. Other known
tensor formats, such as: Tucker, and hierarchical Tucker (HT)
could be also applied \citep{hackbusch2012tensor, khoromskaia2018tensor, khorBook18,  ModelReductionBook15}. 
For the sake of simplicity of exposition, the CP tensor format is considered
in the main part of the paper for illustrative purposes, when direct reference to a low-rank format is
necessary.  The computations are also possible in other compressed formats. 
In \citep{Matthies2020} the computations
of the operations of the algebra were already given for various tensor formats, and they
may also be found e.g.\ in \citep{hackbusch2012tensor, khoromskaia2018tensor, khorBook18, 
 ModelReductionBook15}, so that
we may be brief here.  For numerical experiments we use the Tensor Train 
(TT) software library TT-tool \citep{oseledets2011}.  


In this work we would like to explore the usefulness and applicability of these 
techniques in probability. In probability theory, given a $d$-dimensional RV, 
it is well known that its characterisation is provided by the \pcf \citep{Lukacs70b} (Feller and L\'evi theorems).  Therefore, our contribution will also focus 
on approximation techniques for representing the \pcf in low-rank 
tensor format.  Previous approaches to approximate the \pcf were 
not able to handle high-dimensional RVs, e.g.\ 
\citep{Lukacs70b, Shephard, Witkovsk2016NumericalIO}.  
Our approach can overcome these difficulties.
Other numerical methods which rely on the use of the \pcf may also benefit from the 
proposed approach; e.g.\ the idea to use Markov Chain Monte Carlo (MCMC) methods in the 
Fourier domain to sample from a density proportional to the absolute value of the 
underlying characteristic function is presented in \citep{belomIosipoi2021}.

In a very recent paper
\cite{Dolgov20}, the authors estimate tensor train (TT) ranks for approximated multivariate 
Gaussian \pdfs.  In \cite{DolgovFox19}, the authors use the TT-format to approximate 
multivariate probability distributions.  There they analyse properties of the obtained 
low-rank approximation and use it as a prior distribution in the MCMC 
approach.  In \citep{Witkovsk2016NumericalIO, witkovsky2017computing}, the author 
develops a numerical inversion of characteristic functions. 

Low-rank tensor techniques proved to be very successful in such areas as numerical 
mathematics \citep{Matthies2020, ESHALIMA11_1, Kressner11_low-rank, 
khorBook18, hackbusch2012tensor,  ModelReductionBook15}, 
computational chemistry \cite{khoromskaia2018tensor}, statistics \citep{litv17Tensor} 
and others \citep{GrasedyckKressnerTobler2013}.  Many known types of Green's functions were 
approximated in the low-rank tensor format \citep{khorBook18, khor-lectures-2010, 
Khoromskij_Low_Tacker}, which resulted in drastically reduced computational cost 
and storage.  

General tensor formats and their low-rank approximations, in quantum physics
also known as \emph{tensor networks}, are described in \citep{VI03, Sachdev2010-a, 
EvenblyVidal2011-a, Orus2014-a, BridgemanChubb2017-a, BiamonteBergholm2017}.  
For the mathematical and numerical point of view we refer to 
the review \citep{Kolda:07}, the monographs \citep{HA12, khoromskaia2018tensor, khorBook18},
and to the literature survey on low-rank approximations \citep{GrasedyckKressnerTobler2013}.


%
%
%

Some of the necessary theory is reviewed in \refS{theory}. 
%
The discrete versions of the \pdf and the \pcf are represented in a compressed approximation,
we propose to use low-rank tensor approximations.
The discretised \pdf and \pcf are viewed as tensors equipped with a discrete
version of the point-wise product, the so-called Hadamard product.  With this product
they become (commutative) C$^*$-algebras, so it is possible to define algebraically
the desired functions on these tensors, and on their compressed approximation,
through the use of matrix algebra algorithms \citep{NHigham}. 

Assuming that the \pdf has been represented in such an algebraic setting,
in \refS{statistics} we explain how to compute statistical moments and divergences in 
this discretised framework, using only  the abstract discrete operations, independent
of any particular representation. 
In \refS{tensor-rep} we show as an example how
the algebraic operations may be actually implemented on a low-rank tensor format, namely first
for illustrative purposes on
the simpler canonical polyadic (CP) tensor decomposition, where we 
recall all required definitions and properties.  This tensor format 
is the easiest one to explain the ideas.  
In the example computations later in \refS{Numerics} we actually use the tensor-train
(TT) format, but here the explanation is more complicated.
The actual algorithms---developed for matrices---which can be used
to compute various functions of \pdfs are either iterations or
truncated series expansions, and are listed in \refS{algs}. 
Again,  only  the abstract discrete operations are needed, and the algorithms are independent
of any particular representation.
%
In \refS{Numerics} some numerical examples 
are given, and \refS{concl} concludes.

\section{Theoretical background} \label{S:theory}
The two entities we will be working with are the probability density function (\pdf)
and the corresponding probability characteristic function (\pcf).  These are connected
by the Fourier transform.  The ultimate goal is to use discretised versions of these
two descriptors of RVs to compute for example the \pdf of a mixture model
or the \pdf of the sum of two independent RVs, or quantities of interest (QoIs)
like moments, or other statistics like the relative entropy of a RV, or---for two 
\pdfs---their Kullback-Leibler (KL) or more general $f$-divergences.  While introducing
the notation, this section serves to collect the basic properties of the
different descriptors such as moments,
cumulative distribution function \cdf,  \pdf,  \pcf, 
second characteristic or cumulant generating function \citep{Zoltan13, 
Lukacs70b, Witkovsk2016NumericalIO}, and the moment generating function.
The discretised versions of these QoIs, properties and operations will later be used in the
numerical computations.
In order to evaluate 
point-wise functions of the \pdf, which is necessary to compute various QoIs, 
one has to view the \pdf as well
as the \pcf as elements of function algebras, something which has to be conserved by the
discretisation and low-rank representation.  

\subsection{Notation and summary of basic descriptors}
\label{SS:basic-descr-RV}
%

The random vector
 $  \vX : \Omega\to\C{V}=\Rd $ as
a measurable mapping w.r.t.\ the Borel-$\sigma$-algebra $\F{B}(\Rd)$ on
$\Rd$, defined on a probability space $(\Omega,\F{A},\D{P})$,
where  $\F{A}$ is a $\sigma$-algebra and $\D{P}$ a probability measure was already
introduced in \refS{intro}.  
We are mainly interested in the case where  $\C{V}=\Rd$
is a high-dimensional vector space. 
The associated expectation operator 
is denoted by $\EXP{\cdot}$, the mean as $\vbar{\xi} := \EXP{\vX} = \int_{\Omega} \vX(\omega)
\, \D{P}(\di \omega) \in\Rd$,
and  the mean zero random part as  $\vtil{\xi} := \vX - \vbar{\xi}$.
The canonical Euclidean inner product
 $\bkt{\bx}{\by}_{\C{V}} := \sum_{k=1}^d x_k y_k$ on $\C{V}=\Rd$ is used
to identify the space $\C{V}=\Rd$ with its dual.
The usual order structure on $\Rd=\C{V}$ is assumed, so that
$\bx \le \by$ iff $x_k \le y_k$ for $k=1,\dots,d$, resp.\ $\by-\bx \in \Rd_+=\C{V}_+$,
with positive cone $\Rd_+ = \C{V}_+ = \{\bx\in\Rd \mid x_1\ge 0,\dots, x_d\ge 0 \}$.


The \emph{moments} 
 $\tnb{X}_k$  and the  \emph{central moments} $\tnb{\Xi}_k$ of $\vX$
of order $k\in\N_0$ ---assuming they exist---are basic descriptors of a RV and are denoted as
\begin{equation}   \label{eq:moments-xi}
   \tnb{X}_k = \EXP{\vX^{\otimes k}} \in (\Rd)^{\ten k},\quad
   \tnb{\Xi}_k = \EXP{\vtil{\xi}^{\otimes k}} \in (\Rd)^{\ten k}.
\end{equation}
where for $\bx\in\Rd$ one sets $\bx^{\otimes k} = \Ten_{j=1}^k \bx$.
The second central moment---the \emph{covariance} matrix---is also denoted as
 $\bSigma_{\vX}=\cov{\vX}=\tnb{\Xi}_2 =\tnb{X}_2 - \vbar{\xi}\otimes\vbar{\xi}
 \in \C{V}\ten\C{V}=(\Rd)^{\ten 2}$.
If $\etab$ is another random vector with values in $\C{U}=\Rn$, the
\emph{mixed} and \emph{mixed central} moments are denoted by
\begin{equation}   \label{eq:moments-xi-zeta}
   \tnb{Y}_{k,\ell} = \EXP{\vX^{\ten k}\ten \etab^{\ten \ell}} \quad \text{ and } \quad
   \tnb{\Upsilon}_{k,\ell} = \EXP{\vtil{\xi}^{\ten k}\ten \vtil{\eta}^{\ten \ell}}\in 
       (\Rd)^{\ten k}\ten(\Rn)^{\ten \ell}.
\end{equation}
The \emph{covariance} is also denoted as $\cov(\vX,\etab) = \tnb{\Upsilon}_{1,1}=
\tnb{Y}_{1,1}-\vbar{\xi}\ten\vbar{\eta} \in \C{V}\ten\C{U} = \Rd\ten\Rn$. 

%
The distribution measure
$P_{\vX}$ of $\vX$ for a Borel subset $\C{E}\subseteq\Rd$ 
is as usual  the \emph{push-forward} $\vX_*\D{P}$
of the original measure $\D{P}$:
$ 
   P_{\vX}(\C{E}) := \vX_*\D{P}(\C{E}) := \D{P}(\vX^{-1}(\C{E})) = \EXP{\bbbone_{\C{E}}(\vX)},
$ 
where the \emph{characteristic indicator function} $\bbbone_{\C{E}}(\by)$ 
is unity if $\by\in\C{E}$ and
vanishes otherwise, and is assumed absolutely continuous w.r.t.\ Lebesgue measure on $\Rd$.
This leads to an absolutely continuous 
\emph{cumulative distribution function} (\cdf)
of $\vX$ ---defined for $\by\in\Rd=\C{V}$ via semi-infinite intervals 
$\C{E}_{\by} = \bigtimes_{k=1}^d\; ] \!-\!\infty, y_k ]\in\F{B}(\Rd)$ as
\begin{equation}  \label{eq:distr-fct}
   F_{\vX}(\by): = P_{\vX}(\C{E}_{\by})  = \D{P}(\vX \le \by) = \EXP{\bbbone_{\C{E}_{\by}}(\vX)}.
\end{equation}
The well known properties of $F_{\vX}$ such as positivity $F_{\vX}(\by)\ge 0$,
and the monotonicity
$F_{\vX}(\by_1) \le F_{\vX}(\by_2)$ for $\by_1 \le \by_2$, as well as $F_{\vX}(\by) \to 0$ as
$\by \to -\vek{\infty}$ and $F_{\vX}(\by) \to 1$ as $\by \to +\vek{\infty}$, should
be replicated in any discretised setting.  As $P_{\vX}$ is absolutely continuous,
it has a Radon-Nikod\'ym derivative $p_{\vX}(\by)= \di P_{\vX}/\di \by \in\Lp_1(\Rd,\RR)$ 
w.r.t.\ Lebesgue measure, the \emph{probability density function} (\pdf):
\begin{equation}  \label{eq:prob-dens-fct}
   p_{\vX}(\by)=  \frac{\di}{\di \by}P_{\vX}(\by) = 
    \frac{\dd^d}{\dd y_1\dots\dd y_d} F_{\vX}(\by); \; \text{ with }
    p_{\vX}\ge 0,\;  \text{ and } \int_{\Rd} p_{\vX}(\bx)\, \di \bx = 1.
\end{equation}
These defining relations for the \pdf are directly implied
by the properties of the \cdf $F_{\vX}$, and are important 
to be preserved under discretisation and compressed low-rank approximation.
The positivity  in \refeq{eq:prob-dens-fct} means that geometrically 
speaking densities are in the positive convex
cone of $\Lp_1(\Rd,\RR)$, and the integral relation means that they lie on a 
hyperplane.  The intersection of these two closed convex sets is  denoted by
\begin{equation}  \label{eq:set-prob-dens-fct}
\F{D} := \{ p \mid p \ge 0\} \cap \{ p \mid \int_{\Rd} p(\bx)\, \di \bx = 1\}
\subset \Lp_1(\Rd,\RR), \quad \text{ closed and convex.}
\end{equation}
In addition one notes 
that $\nd{p_{\vX}}_{1} := \int_{\Rd} 
\ns{p_{\vX}(\bx)}\, \di\bx = 1$, i.e.\ $p_{\vX}$ is on the unit sphere.
Convexity of $\F{D}$ in \refeq{eq:set-prob-dens-fct} 
means 
that convex 
combinations of densities $p_{\vX_1}, \dots, p_{\vX_m} \in \F{D}$ 
are again densities; corresponding to \emph{mixture models}. 

Direct \emph{quantities of interest} (QoIs) are usually expected values of functions of
$\vX$, i.e.\ quantities like $\tns{g} = \EXP{\vek{g}(\vX)} \in \Rm=\C{Y}$. 
Obviously the moments  and central moments \refeq{eq:moments-xi}
are special cases  of this by taking as QoI $\vek{g}_k:\bx\mapsto\bx^{\otimes k}$, 
as is the characteristic function \refeq{eq:def-char-fct}.
These quantities may alternatively be computed by integrating over the \pdf $p_{\vX}$ and $\Rd$,
like \refeq{eq:exmpl-diff-entrop} in \refS{intro}:
\begin{equation}   \label{eq:P-qoi}
   \tns{g} = \EXP{\vek{g}(\vX)}= 
   \int_{\Rd} \vek{g}(\bx) \di F_{\vX}(\bx) 
   = \int_{\Rd} \vek{g}(\bx) p_{\vX}(\bx) \di \bx =: \D{E}_{p_{\vX}}(\vek{g}) \in \C{Y}.
\end{equation}

As is well known (e.g.\ \citep{Segal1978}), the Banach space $\Lp_1(\Rd,\RR)$ is
a commutative Banach algebra when the space is equipped with the convolution product:
\begin{equation}  \label{eq:conv-dens}
  (p*q)(\by) := \int_{\Rd} p(\by-\bx)q(\bx)\, \di \bx, \quad \text{ for } p, q \in \Lp_1(\Rd,\RR).
\end{equation}
Observe that if $p$ and $q$ are density functions, so is $p*q$.  This
reflects the well known fact that if $p_{\vX}$ is the density of the RV $\vX$, and
$p_{\etab}$ is the density of the independent RV $\etab$, then $p_{\vX}*p_{\etab}$ is the 
density of the RV $\vX + \etab$.  It also means that the closed convex set $\F{D}$ in
\refeq{eq:set-prob-dens-fct} is stable under convolution.

\ignore{The last topic to touch on to collect the requirements for the
discretised low-rank representation of the \pdf is to look at what is required for
the computation of e.g.\ the \emph{differential entropy} of a RV.  
Recall that the differential entropy
of a \pdf $p_{\vX}$  is defined as $h(p_{\vX}) := \D{E}_{p_{\vX}}(-\log(p_{\vX})) = 
 - \int_{\Rd} \log(p_{\vX}(\bx)) p_{\vX}(\bx) \, \di \bx$.  Hence, in order to compute
the differential entropy $h(p_{\vX})$, one has to be able to compute \emph{point-wise}
functions of the \pdf; in this case the logarithm.  This and other point-wise function needed
for the various divergences, like the square root, pose quite a challenge for a 
compressed low-rank representation, and will be addressed later in \refS{algs}.

To summarise, for the discrete compressed low-rank representation of a \pdf,
the vector space structure is needed to compute convex combinations and the convolution
algebra structure is needed to compute convolutions, as well as the ability to compute 
point-wise functions like the $\log$, the square root, etc.  Additionally, positivity has to
be checked, as well as the condition that the \pdf integrates to unity.  The last
requirement, the ability to compute the Fourier transform, will be explained in the following.}

%
The \emph{characteristic function} (\pcf) of the RV $\vX$ stated already in \refeq{eq:def-char-fct}
may be seen defined  as a QoI via \refeq{eq:P-qoi} with
$g_{\bt}(\bx) = \exp(\ii \bkt{\bt}{\bx})$, i.e.\ the probabilist's Fourier transform---in 
other fields this is considered as the
non-unitary version of the \emph{inverse} Fourier transform \citep{bracewell} ---of the \pdf.
Well known \citep{Lukacs70b, Witkovsk2016NumericalIO} facts about characteristic 
functions $\vphi_{\vX}:\Rd\to\bbC$
are that they are bounded and uniformly continuous--- $\vphi_{\vX} \in \Ck_{bu}(\Rd,\bbC)$ 
---and satisfy $\ns{\vphi_{\vX}(\bt)} \le \vphi_{\vX}(0)=1$
for all $\bt\in\Rd$, i.e.\ $\nd{\vphi_{\vX}}_\infty =1$.
These conditions mean that \pcfs lie on a hyperplane in the
vector space $\Ck_{bu}(\Rd,\bbC)$, obviously reflecting the hyperplane condition---unit 
integral---for the \pdf.

The \pcf is a complex-valued function, but as 
the Fourier transform of
the real and positive \pdf, it has to satisfy some further constraints.
To be able to properly formulate this, one
defines an anti-linear involution or $\star$-operation ``*'' as
$\vphi_{\vX}^{*}(\bt):=  \bar{\vphi}_{\vX}(-\bt)$,
where the overbar denotes the complex conjugate.  Now, as the \pdf
$p_{\vX}$ in \refeq{eq:def-char-fct} is a real function, this implies 
\citep{bracewell, Lukacs70b} that $\vphi_{\vX}$
is Hermitean, i.e.\ invariant w.r.t.\ the conjugation induced by the *-involution; 
so it satisfies $\vphi_{\vX}^*(\bt) =  \bar{\vphi}_{\vX}(-\bt)=  \vphi_{\vX}(\bt)$.
This real subspace of Hermitean functions will be denoted by
$\C{H}\subset \Ck_{bu}(\Rd,\bbC)$.

As is well known \citep{Segal1978}, the Banach space $\Ck_{b}(\Rd,\bbC)$
together with point-wise multiplication and the $\star$-operation is a commutative C*-algebra,
and the real subspace $\C{H}$ of Hermitean functions is a \emph{real}
sub-algebra.
As $p_{\vX}$ is non-negative \citep{Lukacs70b}, the \pcf
$\vphi_{\vX}$ has to be a positive definite function: for any $n\in\D{N}$ and
distinct points $\{\bt_k\in\Rd\}_{k=1}^n$, the matrix
$\vek{\Phi}_{\vX} = (\vphi_{\vX}(\bt_i - \bt_j))_{i,j}\in\bbC^{n\times n}$ 
is Hermitean positive semi-definite in $\bbC^n$.  The positive definite
functions form a convex cone in $\C{H}$, reflecting the positivity of the \pdf.  
These conditions will not be so easy to directly ascertain on a discretised
and compressed low-rank representation of $\vphi_{\vX}$, and are simpler checked on the 
\pdf $p_{\vX}$.  All this means that the \pcf is in the intersection of
this cone $\C{H}$ and the hyperplane of functions with unit value at the origin in $\C{H}$,
 again a closed and convex set, denoted by 
\begin{equation}  \label{eq:set-char-fct}
\F{C} = \C{H} \cap \{\vphi \mid \vphi(0)=1\}
\cap \{ \vphi \mid \vphi \text{ is positive semi-definite } \} \subset \Ck_{bu}(\Rd,\bbC).
\end{equation}
Thus, real convex combinations of \pcfs
are again \pcf, namely of mixture models,
reflecting the analogous property for \pdfs in $\F{D}$ in
\refeq{eq:set-prob-dens-fct}.

It is also easily seen that the point-wise product of two \pcfs
is again a \pcf, and hence the set $\F{C}$ in \refeq{eq:set-char-fct} is stable under
point-wise products.  This reflects the fact that the point-wise product $\vphi_{\vX} \cdot
\vphi_{\etab}$ of the \pcfs of two independent RVs $\vX$ and $\etab$ is just the \pcf of their sum.
This is the analogue to the previous statement about the convolution of two \pdfs,  
showing that $\Fd(\F{D}) \subseteq \F{C}$.

\subsection{Relations through the Fourier transform}
\label{SS:alg-struc-FT}


%
The Fourier transform is known \citep{bracewell, Lukacs70b} to connect the 
algebraic convolution structure on probability densities described previously
with the corresponding multiplication structure on characteristic functions.
Let $p, q \in \Lp_1(\Rd,\RR)$ be integrable real-valued functions in the
convolution Banach-algebra, and denote their Fourier transforms (FT in \refeq{eq:def-char-fct})
by $\phi=\Fd(p), \psi=\Fd(q) \in \C{H}$ (as $p,q$ are real, their FT is Hermitean), 
then the Fourier transform and its inverse $\iFd$ has the following 
well known  \citep{bracewell, Lukacs70b} property:
\begin{equation}  \label{eq:FT-conv-mult}
  \Fd(p * q) = \Fd(p) \cdot \Fd(q) = \phi \cdot \psi \quad\Leftrightarrow\quad
   p*q = \iFd(\phi \cdot \psi).
\end{equation}
Thus the Fourier
transform (FT) is a real algebra homomorphism between the real convolution Banach-algebra
$\Lp_1(\Rd,\RR)$ and the real multiplication algebra $\C{H}$ of Hermitean functions.

In our context, further well known  \citep{bracewell, Segal1978} 
properties of the FT, which will be needed later, are
embodied in the statement $\Fd(\F{D}) \subseteq \F{C}$ at the end of \refSS{basic-descr-RV},
can be summarised as follows: 
there is a correspondence between real-valued functions in $\Lp_1(\Rd,\RR)$ and
Hermitean Fourier transforms $\C{H}$, and 
the positive convex cone in $\Lp_1(\Rd,\RR)$ corresponds to the convex cone of 
positive-definite functions in $\C{H}$, 
while
the hyperplane with unit integral in $\Lp_1(\Rd,\RR)$ corresponds with the hyperplane with
unit value at the origin in $\C{H}$.  These relations will be used  in the discretised
low-rank setting in the following \refSS{low-rank-RV-th} to ascertain the correctness 
of the approximations.  This means that as the \pdf satisfies $p_{\vX}\in\F{D}$ in
\refeq{eq:set-prob-dens-fct}, so the discretised version will have to satisfy a similar
discrete constraint, and dually, as the \pcf satisfies $\vphi_{\vX}\in\F{C}$ in
\refeq{eq:set-char-fct}, the discretised version of this
quantity will have to satisfy a similar discrete constraint.

Another well known property \citep{bracewell, Segal1978} of the Fourier transform that
will be needed 
is how it connects derivatives and multiplication by the co-ordinates.  We assume that all derivatives appearing in the sequel exist and are well defined.
In fact, from \refeq{eq:def-char-fct} one gleans that
 $   (-\ii\dd_{t_k})\, \vphi_{\vX}(\bt) = 
   \int_{\Rd} x_k \exp(\ii \bkt{\bt}{\bx})p_{\vX}(\bx)\, \di \bx =
   \Fd\left( x_k p_{\vX}(\bx)\right)(\bt) $.
Further, denoting the tensor of $k$-th derivatives by
$ 
\Di^k \vphi_{\vX}(\bt) = \left(\frac{\dd^k}{\dd_{t_{i_1}}\dots\dd_{t_{i_k}}}
   \vphi_{\vX}(\bt)\right), 
$ 
one obtains the well-known relation
\begin{equation} 
   \label{eq:der-char-mom-k}
   (-\ii)^k \Di^k  \vphi_{\vX}(\vek{0}) = 
   \int_{\Rd} \bx^{\otimes k}\, p_{\vX}(\bx)\, \di \bx =
   \Fd\left( \bx^{\otimes k}\, p_{\vX}(\bx)\right)(\vek{0}) = \tnb{X}_k,\quad k\in\D{N}_0.
\end{equation}
Similar relations as \refeq{eq:der-char-mom-k} can be obtained by other
characterising functions.
The \emph{second characteristic} function \citep{Lukacs70b} ---sometimes
also labeled as cumulant generating function (cf.\ \refeq{eq:def-cum-g-fct}) ---whose
derivative tensors of order $k$ are essentially the \emph{cumulants} $\tnb{K}_k$
  of $\vX$, is defined as the point-wise logarithm of the \pcf:
\begin{equation}  \label{eq:def-2-char-fct}
   \chi_{\vX}(\bt):=\log(\vphi_{\vX}(\bt)) = \log\left(\EXP{\exp(\ii\bkt{\bt}{\vX})}\right), \,
   \text{ with } (-\ii)^k\Di^k  \chi_{\vX}(\vek{0}) =:  \tnb{K}_k,\quad k\in\D{N}_0. 
\end{equation}
The relations \refeq{eq:der-char-mom-k} and \refeq{eq:def-2-char-fct} involve
the slightly annoying imaginary unit.  To stay with completely real functions
one may switch from the Fourier transform in \refeq{eq:def-char-fct} or
\refeq{eq:def-2-char-fct} to the Laplace transform, at the price of working
with functions which may not be defined for all $\bt\in\Rd$, but maybe only
in a small neighbourhood around $\vek{0}\in\Rd$.
The \emph{moment generating} function is defined as \citep{Lukacs70b} essentially the
reflected \emph{Laplace transform} of the density, or as evaluation of the
\pcf \refeq{eq:def-char-fct} for purely imaginary arguments:
\begin{equation}  \label{eq:def-mom-g-fct}
   M_{\vX}(\bt):=\EXP{\exp(\bkt{\bt}{\vX})} =
   \int_{\Rd} \exp(\bkt{\bt}{\bx})p_{\vX}(\bx) \, \di \bx = \C{L}_d(p_{\vX})(-\bt) =
    \vphi_{\vX}(-\ii\, \bt),
\end{equation}
where $\C{L}_d(p_{\vX})(\bt) = \int \exp(\bkt{-\bt}{\bx})p_{\vX}(\bx) \, \di \bx$
is the two-sided $d$-dimensional \emph{Laplace} transform of $p_{\vX}$. As in 
\refeq{eq:der-char-mom-k}, one obtains
\begin{equation} 
   \label{eq:der-mom-g-k}
   \Di^k  M_{\vX}(\vek{0}) = 
   \int_{\Rd} \bx^{\otimes k}\, p_{\vX}(\bx)\, \di \bx  = \tnb{X}_k,\quad k\in\D{N}_0.
\end{equation}
Closely related is the \emph{cumulant generating} function \citep{Lukacs70b}, the
point-wise logarithm of the moment generating function $M_{\vX}$ in \refeq{eq:def-mom-g-fct}:
\begin{equation}  \label{eq:def-cum-g-fct}
   K_{\vX}(\bt):=\log(M_{\vX}(\bt)) = \log\left(\EXP{\exp(\bkt{\bt}{\vX})}\right), \,
   \text{ with } \Di^k  K_{\vX}(\vek{0}) =  \tnb{K}_k,\quad k\in\D{N}_0.
\end{equation}


\subsection{Grid functions as tensors and 
      algebras on them} 
\label{SS:low-rank-RV-th}
It was already pointed out in \refS{intro} that we want to discretise both the
\pdf and the \pcf, by representing them on a discrete and finite grid.  
We start with the \pdf.  The first thing usually to do is to centre everything around
the mean $\vbar{\xi}$ of the RV $\vX$, i.e.\ to shift co-ordinates on
$\Rd$ by $\bx \mapsto \bx - \vbar{\xi}$.  Another way of saying this is to state that we
work with the centred RV $\vtil{\xi}$.  We assume from now implicitly that this has been
done.  The values of the RV $\vX$ and hence the support 
of its \pdf will thus be around the origin $\vek{0}\in\Rd$.

%
The fully discrete representation of the \pdf and the \pcf is based on
equi-distant grid vectors 
$\hat{x}_{i_\nu,\nu} = \hat{x}_{1,\nu} + (i_\nu-1) \Delta_{x_\nu}$ (with increment $\Delta_{x_\nu}$)
of size $M_\nu$ 
in each dimension $1\le \nu \le d$ of $\Rd$, which  were
introduced in \refSS{motiv}.
The whole grid was denoted by
$\That{X} = \bigtimes_{\nu=1}^d \vhat{x}_{\nu}$,
and we assume that it covers the support of $p_{\vX}(\bx)$. Observe that
functions are implicitly assumed to be periodic \citep{bracewell} when using the FFT, 
hence the total $d$-dimensional volume covered is
   $ V= \prod_{\nu=1}^d M_\nu \Delta_{x_\nu}$,
and the integration rule is implicitly the iterated trapezoidal rule on a periodic
grid, so each point carries the same integration weight $\frk{V}{N}$.  

As already defined in \refeq{eq:pdf-tensor-1} in \refS{intro},
the notation $\tnb{P}:=p_{\vX}(\That{X})$ denotes the tensor $\tnb{P}\in 
\bigotimes_{\nu=1}^d \R^{M_{\nu}} =: \C{T}$, a finite dimensional space with 
$\dim\C{T} = \prod_{\nu=1}^d M_{\nu} =: N$, the components of which are
the evaluation of the \pdf $p_{\vX}$ on the grid $\That{X}$, and
similarly for any other scalar function $f(\bx)$ on $\Rd$.  
Similarly, any vector
valued function $\vek{g}:\C{V}\to\R^m=\C{Y}$, when evaluated on the grid $\That{X}$,
may be viewed as a tensor $\vek{g}(\That{X}) \in \C{Y}\otimes\C{T}$, and
the grid itself $\That{X}$ can be seen as an $\Rd$ valued
functions evaluated on the grid.

To achieve a similar discrete representation of the \pcf, 
it was mentioned in \refSS{motiv} that one makes use of \refeq{eq:grid-char},
$\tnb{\Phi}= \Fd(\tnb{P})$, a discrete version of the  relation \refeq{eq:def-char-fct}, 
denoting the discrete $d$-dimensional Fourier transform on the 
grid $\That{X}$ again by $\Fd$ for the sake of simplicity.
For the dual grid $\That{T} = \bigtimes_{\nu=1}^d \vhat{t}_{\nu}$ from \refSS{motiv},
 the dual grid vector in each dimension is 
 $\vhat{t}_{\nu} := (\hat{t}_{1,\nu},\dots,\hat{t}_{M_\nu,\nu})$.  
 As is well known \citep{bracewell}, if in dimension $\nu$ one has $L_\nu = M_\nu \Delta_{x_\nu}$ as
\emph{period length} for the equi-distantly spaced primal grid with grid spacing $\Delta_{x_\nu}$,
then $\hat{t}_{M_\nu,\nu} = \uppi/\Delta_{x_\nu}$ is the highest (Nyquist) \emph{frequency},
and the equi-distant spacing of the dual grid in dimension $\nu$ is $2 \uppi / L_\nu$.
It is assumed
that the origin $\vek{0}\in\Rd$ is in the dual grid, $\vek{0}\in\That{T}$, and we denote
the index of the origin in the grid with $\vek{j}^0 = (j_1^0,\dots,j_d^0)$, i.e.\ 
$(\hat{t}_{j_1^0,1},\dots,\hat{t}_{j_d^0,d}) = \vek{0}=(0,\dots,0)$.  
This point will be important in some of the statistics resp.\ QoIs to be
described later in \refS{statistics}.  As before, 
the whole grid can be seen as an order $(d+1)$ tensor 
$\That{T} \in \Rd\otimes\C{T} = \C{V}\otimes\C{T}$.
The \pcf on the dual grid is represented as in  \refeq{eq:grid-char} through the tensor 
$\tnb{\Phi}:=\phi_{\vX}(\That{T})\in\C{T}$.
It is on these grid representations $\tnb{P}:=p_{\vX}(\That{X})$ of the \pdf and
$\tnb{\Phi}:=\phi_{\vX}(\That{T})$ of the \pcf  that we propose to operate on to
compute the desired quantities of interest, e.g.\ the differential entropy
\refeq{eq:exmpl-diff-entrop} in \refS{intro}.
A more comprehensive list of such quantities of interest is given in \refS{statistics}.

It was already noted that even for modest values of dimension $d$ and number of discretisation
points $n$ the total amount of data $N=n^d$ may become huge or even non-manageable, and one
has to resort to some kind of compressed representation resp.\ approximation.  Here we
advocate for low-rank tensor representations to allow for efficient computation,
which will be treated in \refS{tensor-rep}.    In such a representation it becomes difficult
to compute point-wise functions (like the log) of values of particular tensors
required for particular statistics (cf.\ \refS{statistics}), as they are not accessible directly. 
To still be able to do these computations, we will rely on certain algebraic properties,
cf.\ \refS{algs}.  These algebraic relations have been pointed out for the non-discrete
entities in the preceding \refSS{alg-struc-FT}, as well as the rôle of the Fourier transform in 
connecting them.   It is now important to ascertain that such algebraic relations
also hold for the grid-discrete quantities. 

%
The set of tensors for the representations of \pdf and \pcf, 
$\C{T} = \bigotimes_{\nu=1}^d \R^{M_{\nu}} \cong \R^N$ is
clearly a vector space.  The \refS{tensor-rep} will explain how the vector
operations can be numerically performed in the low-rank representation.  This assures
that convex linear combinations of different \pdfs or \pcfs can be computed.  
Furthermore, one makes $\C{T}$ into a Euclidean space \citep{hackbusch2012tensor}
by extending the canonical inner products on the $\R^{M_\nu}$ onto $\C{T}$; i.e.\  for elementary tensors 
$\tnb{r}=\bigotimes_{\nu=1}^d \vek{r}_\nu, \tnb{s}=\bigotimes_{\nu=1}^d \vek{s}_\nu, \in \C{T}$,
it is simply 
 $ \bkt{\tnb{r}}{\tnb{s}}_{\C{T}} := \prod_{\nu=1}^d \bkt{\vek{r}_\nu}{\vek{s}_\nu}_{\R^{M_\nu}} $,
and then extended to all of $\C{T}$ by linearity.

It is convenient to extend this inner product to larger tensor products
\citep{hackbusch2012tensor}, and only perform a partial inner product.  For
example, if $\tnb{S} = \vek{y}\otimes\tnb{r}$ is an elementary tensor in
$\C{Y}\otimes\C{T}$, and $\tnb{u}\in\C{T}$, then the partial product,
denoted as before, is
 $ \bkt{\tnb{S}}{\tnb{u}}_{\C{T}} := \bkt{\vek{y}\otimes\tnb{r}}{\tnb{u}}_{\C{T}} 
  = \bkt{\tnb{r}}{\tnb{u}}_{\C{T}}\, \vek{y} \in \C{Y} $,
and then extended to all of $\C{Y}\otimes\C{T}$ by linearity.  It is practically
a contraction over all indices related to $\C{T}$.

The next task is to introduce the algebraic structures.  For 
$\vek{r} = (r_1,\dots,r_{M_\nu}), \vek{s}=(s_1,\dots,s_{M_\nu}) \in \R^{M_\nu}$, 
we want to define the \emph{circular convolution}
$\bz := (z_1,\dots,z_{M_\nu}) = \vek{r} * \vek{s}\in\R^{M_\nu}$ 
component-wise as
 $ z_k := \sum_{\ell=1}^{M_\nu} r_\ell\, s_{m}, \, m = ((k-\ell)\mod {M_\nu}) + 1;\, 1\le k \le {M_\nu} $
 \citep{bracewell, hackbusch2012tensor}.
This is an associative and commutative product; the convolution
algebra on $\R^{M_\nu}$.  It is the discrete version of the convolution algebra $\Lp_1(\Rd,\R)$
considered in \refSS{alg-struc-FT}.  One should point out \citep{bracewell}
that two steps are involved in going from $\Rd$ resp.\ $\R$ to the finite grids $\That{X}$
resp.\ $\vhat{x}_\nu$:  One is the truncation of the infinite domain, say $\R$ in dimension
$\nu$, to a finite one---the interval $[\xi_\nu^{(\min)}, \xi_\nu^{(\max)}]$.  
This picks only certain
discrete frequencies or wavenumbers \citep{bracewell} from the continuum in the
Fourier transform, and they are all multiples of a basic frequency resp.\ wavenumber.
  The second step is the
use of a finite grid, this picks a finite number of the discrete frequencies.
The effect of this is to make everything implicitly periodic, through periodic continuation.
This is reflected in the use of the \emph{circular} convolution.

Having defined the circular convolution on each $\R^{M_\nu}$, it is 
defined on $\C{T} = \bigotimes_{\nu=1}^d \R^{M_{\nu}}$ first via  
elementary tensors $\tnb{r}=\bigotimes_{\nu=1}^d
\vek{r}_\nu,  \tnb{s}=\bigotimes_{\nu=1}^d \vek{s}_\nu\; \in \C{T}$ as 
 $ \tnb{z} = \tnb{r} * \tnb{s} := \bigotimes_{\nu=1}^d \vek{r}_\nu * \vek{s}_\nu $,
and extended to all of $\C{T}$ by linearity.  This again is an associative and
commutative product; the convolution algebra  on $\C{T}$.  In passing one may remark
that the statement that the discrete version $\tnb{\Phi}$ of the \pcf is Hermitean
positive definite is equivalent with the statement that the linear (in $\tnb{r}$) operator
$\tnb{K}_{\tnb{\Phi}}:\tnb{r}\mapsto\tnb{\Phi} * \tnb{r}$ is Hermitean or self-adjoint
and positive definite; $\tnb{K}_{\tnb{\Phi}}$ is the operator of convolution with $\tnb{\Phi}$.

The 
discrete version of the point-wise product 
is the Hadamard product, and again first formulated component-wise on $\R^{M_\nu}$:
For $\vek{r} = (r_1,\dots,r_{M_\nu}), \vek{s}=(s_1,\dots,s_{M_\nu}) \in \R^{M_\nu}$ it is denoted by 
 \citep{hackbusch2012tensor} $\bz := (z_1,\dots,z_{M_\nu}) = \vek{r} \odot \vek{s}\in\R^{M_\nu}$ 
and defined component-wise as 
$  z_k :=  r_k \cdot s_k, \; 1 \le k \le {M_\nu} $.
This also is an associative and commutative product; the Hadamard
algebra  on $\R^{M_\nu}$.  
It is again defined to $\C{T} = \bigotimes_{\nu=1}^d \R^{M_{\nu}}$ via elementary tensors 
$\tnb{r}=\bigotimes_{\nu=1}^d \vek{r}_\nu$ and 
$\tnb{s}=\bigotimes_{\nu=1}^d \vek{s}_\nu \in \C{T}$ as
 $ \tnb{z} = \tnb{r} \odot \tnb{s} := \bigotimes_{\nu=1}^d \vek{r}_\nu \odot \vek{s}_\nu $,
and extended to all of $\C{T}$ by linearity.  And again this is an associative and
commutative product; the Hadamard algebra  on $\C{T}$.  
The interaction with the inner product is quite simple,
it is an elementary calculation to verify that for $\tnb{w}, \tnb{r}, \tnb{s}\in\C{T}$
 $ \bkt{\tnb{w}\odot\tnb{r}}{\tnb{s}}_{\C{T}} =  \bkt{\tnb{r}}{\tnb{w}\odot\tnb{s}}_{\C{T}}$,
which means that the linear (in $\tnb{r}$) operation of Hadamard multiplication by $\tnb{w}$,
i.e.\ $\tnb{L}_{\tnb{w}}:\tnb{r}\mapsto\tnb{w}\odot\tnb{r}$, is self-adjoint.


It is easily seen that the Hadamard algebra has a multiplicative unit, which
we denote by $\tnb{1} = (1_{i_1,\dots,i_d})$ --- the tensor with all ones --- satisfying
$\tnb{r}\odot\tnb{1}=\tnb{r}$ for any $\tnb{r}\in\C{T}$.  Defining
$\vek{1}_\nu := (1,\dots,1) \in \R^{M_\nu}$, a vector of all ones in dimension $\nu$ (which is
the Hadamard unit for the Hadamard algebra on $\R^{M_\nu}$), it is
not difficult to see that the Hadamard unit on the tensor product has rank one: 
$\tnb{1} = \bigotimes_{\nu=1}^d \vek{1}_\nu$; it is the simple tensor product of the 
Hadamard units on each $\R^{M_\nu}$.

The unit further allows to introduce a discrete expectation or integral operator.
Recall that we have a tensor quadrature grid, equi-spaced in each dimension, in order to use
the FFT.  As mentioned before, each point has the integration weight $\frk{V}{N}$. 
Hence for a tensor
$\tnb{P}$, representing a function $p(\bx)$ evaluated on the grid, the approximate integral is
\begin{equation}  \label{eq:discr-int}
  \int p(\bx)\,\di\bx \approx \C{S}(\tnb{P}) := \frac{V}{N} \bkt{\tnb{P}}{\tnb{1}}_{\C{T}} .
\end{equation}
\refeq{eq:discr-int} is only a convenient way of writing the approximate integral,
and obviously there is no need to actually compute the inner product with $\tnb{1}$, i.e.\ multiply
each entry in the low-rank representation of $\tnb{P}$ with unity.
This \refeq{eq:discr-int} can be simply extended to discrete integrands of the 
form $\tnb{S}\in\C{Y}\otimes\C{T}$ through the use of the partial inner product 
defined above with the result $\C{S}(\tnb{S})\in\C{Y}$.
If $\tnb{F}$ is a tensor which represents the grid-values of a function $f(\bx)$, i.e.\
$\tnb{F} = f(\That{X})$, and $\tnb{P}$ represents the \pdf $p_{\vX}$,
one may define a discrete version of the expectation, which can be 
extended to other tensors $\tnb{S}\in\C{Y}\otimes\C{T}$, by
\begin{equation}  \label{eq:discr-exp}
 \EXP{f(\vX)}=\D{E}_{p_{\vX}}(f) = \int_{\Rd} f(\bx) p_{\vX}(\bx)\, \di \bx \; \approx \;
   \C{S}(\tnb{F}\odot\tnb{P}) =   \frac{V}{N} \bkt{\tnb{F}}{\tnb{P}}_{\C{T}} 
   =: \D{E}_{\tnb{P}}(\tnb{F}).
\end{equation}

%
To translate the statements in 
\refSS{alg-struc-FT} into the present discrete setting, let $\tnb{p}, \tnb{q} \in \C{T}$
be two tensors, and $\tnb{\phi} = \Fd(\tnb{p}), \tnb{\psi}=\Fd(\tnb{q}) \in \C{T}$
their discrete Fourier transforms.  Then one has just as in \refeq{eq:FT-conv-mult}
\begin{equation}  \label{eq:FT-conv-mult-tens}
  \Fd(\tnb{p} * \tnb{q}) = \Fd(\tnb{p}) \odot \Fd(\tnb{q}) = \tnb{\phi} \odot \tnb{\psi}
   \quad\Leftrightarrow\quad \tnb{p}*\tnb{q} = \iFd(\tnb{\phi} \odot \tnb{\psi}),
\end{equation}
showing that the discrete Fourier transform (FT) is an algebra isomorphism between the
convolution algebra $(\C{T},*)$ and the Hadamard algebra $(\C{T},\odot)$.
This makes it relatively easy to compute the discrete convolution of the tensor
representations of two densities, say $p$ represented by $\tnb{P}$, and
$q$ represented by $\tnb{Q}$:  Then the density $p*q$ corresponds to $\tnb{P}*\tnb{Q}$,
computed with their Fourier transforms $\tnb{\Phi}=\Fd(\tnb{P})$ and
$\tnb{\Psi}=\Fd(\tnb{Q})$ as
 $ \tnb{P}*\tnb{Q} = \iFd(\tnb{\Phi}\odot\tnb{\Psi}) = \iFd(\Fd(\tnb{P})\odot\Fd(\tnb{Q}))$.

For a real tensor $\tnb{w}$, the Hadamard multiplication operator $\tnb{L}_{\tnb{w}}$
of Hadamard multiplication with $\tnb{w}$ was defined above.  As the Hadamard
algebra is commutative, the Hadamard multiplication operators commute with each
other.  Thus it is theoretically clear \citep{Segal1978} that they can be
simultaneously diagonalised.  But it is elementary to see this explicitly.
If an arbitrary tensor $\tnb{r}\in\C{T}$ were written as a vector, the 
action $\tnb{L}_{\tnb{w}}(\tnb{r})=\tnb{w}\odot\tnb{r}$ would be the action of a 
diagonal matrix, and diagonal matrices obviously commute.
This means that $\tnb{L}_{\tnb{w}}$ is fully diagonalised, 
and the components $\tns{w}_{j_1,\dots,j_d}$ of
the tensor $\tnb{w}$ are in fact the diagonal elements, i.e.\ the eigenvalues
of $\tnb{L}_{\tnb{w}}$.  Obviously, to an eigenvalue $\tns{w}_{j_1,\dots,j_d}$
belongs the canonical unit vector with the same index as eigenvector, in tensor
notation the eigenvector is $\tnb{v}^{(j_1,\dots,j_d)} = 
(\updelta_{(j_1,\dots,j_d),(i_1,\dots,i_d)})_{i_1,\dots,i_d}$,
i.e.\ $\tnb{L}_{\tnb{w}}(\tnb{v}^{(j_1,\dots,j_d)}) = 
\tns{w}_{j_1,\dots,j_d}\tnb{v}^{(j_1,\dots,j_d)}$.

Let $\tnb{\phi} = \iFd(\tnb{w})$, and let  $\tnb{r}\in\C{T}$ be any tensor.
The convolution operator $\tnb{K}_{\tnb{\phi}}$ of convolution with a
Hermitean $\tnb{\phi}$ was defined above.  One then has 
\begin{align}  \label{eq:conv-EV-expl-1}
  \tnb{K}_{\tnb{\phi}}(\tnb{r}) &= \tnb{\phi}*\tnb{r} = \iFd(\Fd(\tnb{\phi})\odot\Fd(\tnb{r}))
  = \iFd(\tnb{w}\odot\Fd(\tnb{r}))   = 
  \iFd(\tnb{L}_{\tnb{w}}(\Fd(\tnb{r}))),   \quad \text{ or} \\    \label{eq:conv-EV-expl-2}
  \tnb{L}_{\tnb{w}}(\tnb{r}) &= \tnb{w}\odot\tnb{r} = \Fd(\tnb{\phi})\odot\tnb{r}
  = \Fd(\tnb{\phi}*\iFd(\tnb{r}))   = \Fd(\tnb{K}_{\tnb{\phi}}(\iFd(\tnb{r}))).
\end{align}
This shows that $\tnb{K}_{\tnb{\phi}}$ and $\tnb{L}_{\tnb{w}}$ are unitarily
equivalent and have the same eigenvalues; and since those of $\tnb{L}_{\tnb{w}}$
are the components $\tns{w}_{j_1,\dots,j_d}$ of the tensor $\tnb{w}$,
those of $\tnb{K}_{\tnb{\phi}}$ are the same, i.e.\ the components of 
$\tnb{w} = \Fd(\tnb{\phi})$, the Fourier transform of the Hermitean
convolution tensor $\tnb{\phi}$.

When manipulating tensors which represent discrete versions of the \pdf resp.\
the \pcf, one wants to make sure that these manipulations do not destroy the fundamental
properties of these objects.  Let us consider discrete densities first, and assume that
$\tnb{P}\in\C{T}$ is the discrete representation of a density.  Then we should expect
$\tnb{P} \ge \tnb{0}$, i.e.\ for each $1 \le i_\nu \le M_\nu,\; 1\le\nu\le d$: 
$\tns{P}_{i_1,\dots,i_d} \ge 0$.
We have tacitly assumed that $\tnb{P}$ is also real-valued; this is not necessarily so
if it is numerically computed with the involvement of the FT.  But this condition is easy to
check, as usually the real and imaginary parts of a complex tensor are approximated
separately; one only has to make sure the imaginary part is identical to zero---or
even better, not computed or stored at all by using a half-length compressed real
version of the FFT.  This incidentally insures that its FT $\tnb{\Phi} = \Fd(\tnb{P})$
is Hermitean.  This condition is preserved by the discrete FT.
Coming back to the positivity of $\tnb{P}$, this can be checked by
computing its \emph{minimum} component $\tns{P}_{\min}$; this minimum should be
non-negative, making $\tnb{L}_{\tnb{P}}$ positive.  
If it happens to be negative, those components can be removed and
set to zero.  It will be explained later in \refS{algs} on how to do this.  Observing this will also
automatically make its FT $\tnb{\Phi}$ resp.\ $\tnb{K}_{\tnb{\phi}}$ positive definite.

A further condition is that the density integrates to unity.  Hence in the discrete
setting we should expect $\C{S}(\tnb{P}) = 1$.  Incidentally, if $\tnb{\Phi}
=\Fd(\tnb{P})$ is its FT, then 
this is equivalent with $\tnb{\Phi}_{\vek{j}^0} = \tns{\Phi}_{j_1^0,\dots,j_d^0} = 1$.
If $\tnb{P}$ does not meet this condition, then it can be rescaled appropriately.

%
The last task to address 
in connection with compressed resp.\ low-rank tensor approximations
is how to compute point-wise functions $f(\tnb{w})$ of a real tensor $\tnb{w}\in\C{T}$,
which means the same function applied to each component.  
Normally, in a full representation, this computation is no problem,
but in a compressed representation the component values are not directly accessible.
The Hadamard algebra can be used to accomplish this
function evaluation, continuing \citep{ESHALIMA11_1, Matthies2020}.
This builds theoretically on generally well known results (e.g.\ \citep{Segal1978}) 
about the evaluation of functions of self-adjoint linear operators like 
$\tnb{L}_{\tnb{P}}$ with the spectral functional calculus, and abstractly 
on how to compute functions of self-adjoint elements in an Abelian C*-algebra.  
More specifically, as will be seen later, one may use well known 
algorithms for (real and self-adjoint) matrices \citep{NHigham} to actually do the calculations.

The simplest functions like linear combinations follow from the vector space structure of $\C{T}$.
The next in complexity are powers and polynomials.  For $m\in\N_0$ the powers are defined as usual by
$\tnb{r}^{\odot m}:= \tnb{r}^{\odot (m-1)}\odot\tnb{r}$, setting $\tnb{r}^{\odot 0}=\tnb{1}$.
This way one may evaluate \emph{polynomials} $f_p(t) = \sum_{m=0}^M \beta_m t^m$
by replacing $t^m$ with $\tnb{r}^{\odot m}$. 

The inverse function $f_{i}(t) = 1/t = t^{-1}$ is easy to deal with, as
the existence of a unit allows the definition of an inverse element:  $\tnb{r}\in\C{T}$
is called invertible, if there exists a---unique---element $\tnb{w}\in\C{T}$ such
that $\tnb{r}\odot\tnb{w}=\tnb{1}$; it is denoted by $\tnb{w}=\tnb{r}^{\odot-1}$.
It is obvious that for $\tnb{r}\in\C{T}$ to have a Hadamard inverse, no component
 can vanish, $\tns{r}_{i_1,\dots,i_d}\ne 0$. 
and in that case $\tnb{r}^{\odot-1} = (1/\tns{r}_{i_1,\dots,i_d})$.  This is up to
now only a definition of $f_{i}(\tnb{r})=\tnb{r}^{\odot-1}$, and not yet an algebraic 
way of computing it.

Power series $f_{ps}(t) = \sum_{m=0}^\infty \beta_m t^m$,
or more generally $f_{ps}(t)= \sum_{m=0}^\infty \beta_m (t-t_0)^m$ 
are---thanks to the Cayley-Hamilton theorem---actually polynomials on the 
finite dimensional Hadamard algebra---see the remarks on general functions in the following.
In case $f_{h}$ is a holomorphic function in a complex domain containing the
values of $\tnb{r}$,  it can also be evaluated in the algebra via Cauchy's formula:
 $ f_{h}(\tnb{r}) = \oint_{\Gamma} f_{h}(z)\left(z\cdot\tnb{1} - \tnb{r} \right)^{\odot -1}
    \, \di z \in \C{T} $,
where $\Gamma$ is a contour in $\D{C}$ inside the domain of holomorphy,
with all values of $\tnb{r}$ inside the contour.   



The computation of a general function $f(\tnb{r})$ ---where
$f$ is a real valued function defined on a subset of $\R$ which includes all the
values of $\tnb{r}$ ---makes use of
the representation $\tnb{r} \mapsto \tnb{L}_{\tnb{r}}$ from the algebra $(\C{T},\odot)$
into the algebra $(\E{L}(\C{T}),\circ)$ ---the algebra of linear operators with
concatenation ``$\circ$'' as multiplication.  This is an injective algebra homomorphism onto
an Abelian sub-algebra of $\E{L}(\C{T})$ of self-adjoint operators.
The general way to compute a function $f(\tnb{r})$ of a self-adjoint $\tnb{L}_{\tnb{r}}$
is to use the spectral calculus (e.g.\ \citep{Segal1978}).
This will work in any unital C*-algebra; it is mentioned only
for general background orientation, the situation here is much simpler. 
We saw that the operators $\tnb{L}_{\tnb{r}}$ can be represented by diagonal matrices.
This means that \emph{matrix algorithms}, which only use the matrix algebra operations,
can be used to compute $f(\tnb{L}_{\tnb{r}})$ to represent $f(\tnb{r})$, 
see e.g.\ \citep{NHigham}.
The algorithms will explicitly be addressed later in \refS{algs}, and may be used
directly on the Hadamard algebra, and there is no need to actually use matrices.
With that in mind, it may be mentioned in passing that as the
spectrum resp.\ set of eigenvalues of $\tnb{L}_{\tnb{r}}$ is a finite set, any function
$f(\tnb{r})$ could in principle be interpolated by a polynomial.  But this is not practical,
as the polynomial would in general have degree $N=\dim \C{T}$, which we assume to be a huge
number, so that typically the matrix algorithms are more economical.  

\section{Computation of moments and divergences}  \label{S:statistics}
%
{This section explains how to compute QoIs already mentioned in \refS{theory}
from \pdfs or \pcfs given in compressed repesentation.}  
It is thereby assumed that point-wise functions of the compressed representation resp.\
low-rank tensor approximation can be computed.  Algorithms for this task will be shown
in  \refS{algs}.  Moreover, the computation of various
distances and divergences between \pdfs will be treated in this
discrete setting, noting the point-wise functions which are required to compute them.


%
%

%
If it is desired to compute some QoI---the expected value of a function 
$\vek{g}(\vX)$ (with $\vek{g}:\Rd\to\C{Y}$) of the RV $\vX$ like in 
\refeq{eq:P-qoi}, then it is first
necessary to represent this function in low-rank format on the grid $\That{X}$,
i.e.\ to find $\tnb{G} = \vek{g}(\That{X})\in\C{Y}\otimes\C{T}$. 
If $\tnb{P}$ represents the \pdf $p_{\vX}$,
the discrete version  of the expectation  was already given in
\refeq{eq:discr-exp}:
$  \D{E}_{\tnb{P}}(\tnb{G}) =  \frk{V}{N} \bkt{\tnb{G}}{\tnb{P}}_{\C{T}} $.
Moments are a special case: 
for the mean one takes
the tensor $\tnb{x}$, 
or for higher moments the
tensor $\tnb{x}^{\otimes k}$: 
 $ \vbar{\xi} = \tnb{X}_1 \approx  \D{E}_{\tnb{P}}(\tnb{x})$,  and  
 $ \tnb{X}_k \approx \D{E}_{\tnb{P}}(\tnb{x}^{\otimes k})$. 
For mixed moments these relations can be used in an analogous fashion.


Other characterising functions were sketched at the end of \refSS{alg-struc-FT}.
The simplest one is the second characteristic function $\chi_{\vX}(\bt) = \log(\vphi_{\vX}(\bt))$
from \refeq{eq:def-2-char-fct}.  Given the tensor $\tnb{\Phi}$ representing $\vphi_{\vX}$,
this is just the point-wise logarithm 
\begin{equation}  \label{eq:discr-2-char-fct}
        \tnb{\Upsilon} = \log(\tnb{\Phi}), \quad \text{ such that } 
        \tnb{\Upsilon} = \chi_{\vX}(\That{T}).
\end{equation}
From this one may compute cumulants 
\begin{equation}  \label{eq:comp-discr-cum}
  \tnb{K}_k \approx 
  \left(\Fd\left(\tnb{x}^{\otimes k}\odot  \iFd(\tnb{\Upsilon}) \right)\right)_{\vek{j}^0} . 
\end{equation}

The moment generating function $M_{\vX}(\bt)=\EXP{\exp\bkt{\bt}{\vX}}$
\refeq{eq:def-mom-g-fct} can in principle be generated in the same way as
the characteristic function $\vphi_{\vX}(\bt)=\EXP{\exp(\ii \bkt{\bt}{\vX})}$
with the help of the RV $\vX(\cdot)$
to give its tensor representation 
$\tnb{M}=M_{\vX}(\That{T})$ on the grid $\That{T}$.  
Its derivatives are directly the moments, so with the FT one has an
alternative approximate computation of the moments:
\begin{align}  \label{eq:comp-discr-mom-3}
  \vbar{\xi} &\approx  \D{E}_{\C{T}}(\tnb{x}) 
   = \ii \,\left(\Fd\left( \tnb{x}\odot \iFd(\tnb{M})\right)\right)_{\vek{j}^0}, 
  \quad \text{ and } \\  \label{eq:comp-discr-mom-3a}
  \tnb{X}_k &\approx \D{E}_{\C{T}}(\tnb{x}^{\otimes k}) = 
  \ii^k\,\left(\Fd\left( \tnb{x}^{\otimes k}\odot \iFd(\tnb{M})\right)\right)_{\vek{j}^0}. 
\end{align}
Similarly to the discrete version of second characteristic function \refeq{eq:discr-2-char-fct},
the cumulant generating function $K_{\vX}(\bt)=\log(M_{\vX}(\bt))$ in \refeq{eq:def-cum-g-fct}
also has a discrete version, the point-wise logarithm of the discrete moment generating
function $\tnb{M}$: 
\begin{equation}  \label{eq:discr-2-cum-g-fct}
        \tnb{Z} = \log(\tnb{M}), \quad \text{ such that } 
        \tnb{Z} = K_{\vX}(\That{T}).
\end{equation}
Analogous as before in \refeq{eq:comp-discr-mom-3a}, the cumulants may be expressed via
\begin{equation}  \label{eq:comp-discr-cum-2}
  \tnb{K}_k \approx 
  \ii^k\,\left(\Fd\left( \tnb{x}^{\otimes k}\odot \iFd(\tnb{Z})\right)\right)_{\vek{j}^0}. 
\end{equation}
%

If the discretised \pcf $\tnb{\Phi}$ has 
low tensor rank, then  $\tnb{X}_k$ will also be of a low-rank.
Note that we cannot say the same about $\tnb{K}_k$ from \refeq{eq:comp-discr-cum} 
(or \refeq{eq:comp-discr-cum-2}). The tensor rank of $\tnb{K}_k$ depends on the tensor rank 
of $\tnb{\Upsilon}$ (or $\tnb{Z}$), and the latter can have a high rank due to the 
involved $\log(\cdot)$ function.  Hence it is assumed here that the series expansion 
described in \refSS{series} will allow it to approximate $\tnb{\Upsilon}$ 
(or $\tnb{Z}$) in a low-rank tensor format.

%
%
\begin{center}
\begin{table}[htbp!]
 \centering
  \caption{List of some typical divergences and distances.}
 \label{table:divergences}  
 \begin{tabular}{ll}
 \hline
Name of the divergence & $D_{\bullet}(p \Vert q)$ \\ \hline
Kullback–Leibler-(KL) --- $D_{\text{KL}}$: &	$ {\displaystyle \int \left(
        \log(p(\bx)/q(\bx))\right) p(\bx)\, \di \bx = \D{E}_p(\log(p/q))}$\\ 
sqrd.\ Hellinger dist. --- $(D_{H})^2$: &$  {\displaystyle \frac{1}{2}
    \int \left(\sqrt{p(\bx)} - \sqrt{q(\bx)}\right)^2 \, \di \bx  }$ \\
Bregman divergence  --- $D_{\phi}$:	& $ {\displaystyle \int \left[ (\phi(p(\bx)) - \phi(q(\bx)))
         - (p(\bx) - q(\bx)) \phi^{\prime}(q(\bx))\right]\, \di \bx }$ \\
Bhattacharyya distance --- $D_{\text{Bh}}$: & $  {\displaystyle -\log\left(
 \int \sqrt{(p(\bx) q(\bx))} \, \di \bx \right) }$ \\ \hline
\end{tabular}
\end{table}
\end{center}

Let us now turn to quantities which characterise the whole distribution, or
the difference or distance between them.  
A simple characterisation of a distribution is its entropy, which was
already encountered in \refS{intro}.  For a continuous distribution
$p$ this is the \emph{differential entropy}, requiring the 
point-wise logarithm of $\tnb{P}$:
\begin{equation}  \label{eq:discr-diff-entrop}
   h(p) := \D{E}_p(-\log(p)) \approx \D{E}_{\tnb{P}}( -\log(\tnb{P}))
    = - \frac{V}{N} \bkt{\log(\tnb{P})}{\tnb{P}},
\end{equation}
where the expectation $\D{E}_{\tnb{P}}(\cdot)$ is computed as in \refeq{eq:discr-exp}.

\begin{center}
\begin{table}[htbp!]
 \centering
  \caption{Discrete approximations for divergences listed in Table~\ref{table:divergences}. Formulas for computations are given in \refS{theory} and \refS{algs}.}
 \label{table:divergences-discr}  
 \begin{tabular}{ll}
 \hline
Name of the divergence & Approximation of $D_{\bullet}(p \Vert q)$ \\ \hline
Kullback–Leibler-(KL) --- $D_{\text{KL}}$: &	$ {\displaystyle 
  \frac{V}{N} (\bkt{\log(\tnb{P})}{\tnb{P}} - 
              \bkt{\log(\tnb{Q})}{\tnb{P}}) }$\\
squared Hellinger distance --- $(D_{H})^2$: &$  {\displaystyle 
    \frac{V}{2 N} \bkt{\tnb{P}^{\odot 1/2} - \tnb{Q}^{\odot 1/2}}{\tnb{P}^{\odot 1/2} - \tnb{Q}^{\odot 1/2}}}$ \\
Bregman divergence  --- $D_{\phi}$:	& $ {\displaystyle \C{S} \left( 
  (\phi(\tnb{P}) - \phi(\tnb{Q})) - (\tnb{P} - \tnb{Q})\odot \phi^{\prime}(\tnb{Q})\right)}$ \\
Bhattacharyya distance --- $D_{\text{Bh}}$: & $  {\displaystyle   
 -\log \left(\frac{V}{N} \bkt{\tnb{P}^{\odot 1/2}}{\tnb{Q}^{\odot 1/2}}\right)}$ \\ \hline
\end{tabular}
\end{table}
\end{center}

To compare two \pdfs $p$ and $q$, one uses divergences or
distances between them.  The best known is probably the relative entropy, also known as
the Kullback-Leibler (KL) divergence.
Divergences are a kind of generalisation
of distances, but unlike distances they do not have to be symmetric.
Let $p, q \in \F{D}$ be two \pdfs, and $\tnb{P}, \tnb{Q}\in\C{T}$
their low-rank tensor representations.
Some well known divergences and distances  are given
in \reftab{table:divergences}.
For the Bregman divergence, $\phi$ has to be a real convex function, 
e.g.\ $\phi(t)=t^2$.  These divergences and distances may be computed by
the discrete approximations shown in \reftab{table:divergences-discr}.

It is evident that one has to be able to compute the point-wise logarithm and 
the square root of a low-rank tensor.
For the Bregmann divergence one has to be able to compute  point-wise the
convex function $\phi$ and its derivative $\phi^\prime$ of a low-rank tensor
representation, an easy task for $\phi(t)=t^2$.

%
%
The $f$-divergence is a way to define many divergences in a unifying way
\citep{LieseVajda06, NielsenNock2013}, depending on a convex function $f$.
One requires that $f$ be a convex function, 
satisfying  $f(1) = 0$.  Then 
the $f$-divergence of $p$ from $q$ and its discrete approximation is defined as
\begin{equation}  \label{eq:f-diverg}
        D_f ( p \Vert q ) := {\displaystyle \D{E}_{q}\left(f\left( \frac{p}{q} 
        \right) \right) }
        \approx \D{E}_{\tnb{Q}}\left(f(\tnb{P}\odot\tnb{Q}^{\odot -1})\right) 
 =\frac{V}{N} \bkt{f(\tnb{P}\odot\tnb{Q}^{\odot -1})}{\tnb{Q}}.
\end{equation}

Many common divergences, such as the KL-divergence, the Hellinger distance, and the
total variation distance, are special cases of the $f$-divergence, coinciding 
with a particular choice of $f$. In \reftab{table:f-divergence} the functions $f$ for
some common  $f$-divergences \citep{LieseVajda06, NielsenNock2013} are listed.

\begin{center}
\begin{table}[h]
 \centering
  \caption{Some examples  of  the function $f$ for the  $f$-divergence.}
 \label{table:f-divergence}  
 \begin{tabular}{ll}
 \hline
Name of the divergence & Corresponding $f(t)$  \\ \hline
KL-divergence &	$ {\displaystyle t \log(t) }$\\
reverse KL-divergence  &$  {\displaystyle -\log(t)}$ \\
squared Hellinger distance 	& $ {\displaystyle ({\sqrt{t}}-1)^{2}}$ \\
total variation distance 	& $ {\displaystyle \ns{t-1}/2}$ \\
Pearson $\chi^2_P$-divergence	& $ {\displaystyle (t-1)^{2}}$ \\
Neyman $\chi^2_N$-divergence (reverse Pearson) & $ {\displaystyle t^{-1}-1}$ \\
Pearson-Vajda $\chi^k_P$-divergence	& $ {\displaystyle (t-1)^{k}}$ \\
Pearson-Vajda $\ns{\chi}^k_P$-divergence	& $ {\displaystyle \ns{t-1}^{k}}$ \\
Jensen-Shannon-divergence   & $ {\displaystyle t\log(t)-(t+1)\log((t+1)/2) }$ \\ \hline
\end{tabular}
\end{table}
\end{center}

Hence to compute the discrete approximation, one has to be able to compute not
only the Hadamard inverse of a low rank tensor, but also the function $f$.
These are essentially the ones which have been discussed already, the only new
one is the absolute value.

\newcommand{\CP}{\mrm{CP}}
\newcommand{\Tuck}{\mrm{T}}
\newcommand{\TT}{\mrm{TT}}

\section{Tensor formats}      \label{S:tensor-rep}
%
{In this section some well-known technical details scattered in the literature
about algebraic computations  in low-rank tensor formats are collected.} 
Many tensor formats are used in quantum physics under the name \emph{tensor networks},
see \citep{VI03, Sachdev2010-a, EvenblyVidal2011-a, Orus2014-a, 
BridgemanChubb2017-a, BiamonteBergholm2017}. 
The canonical polyadic (CP) \citep{Hitch:27} and \emph{Tucker} \citep{Tuck:66} formats 
have been well known for a long time and are therefore very
popular. For instance, CP and Tucker rank-structured tensor formats have been 
applied in chemometrics and in signal processing \citep{smilde-book-2004, Cichocki:2002}. 
The \emph{Tensor Train} format was originally developed in quantum physics and
chemistry as ``matrix product states'' (MPS), see \citep{VI03} and references therein, 
and rediscovered and developed further in \cite{oseledetsTyrt2010, oseledets2011,
BallaniGrasedyck2015,ds-alscross-2017pre} as tensor train format.
The hierarchical tensor (HT) format was introduced in \citep{HackbuschKuhn2009}, and further 
considered in \citep{Grasedyck2010}.

The CP format is cheap, it is simpler than the Tucker or TT format, but, 
compared to others, there are no reliable algorithms to compute CP decompositions 
for $d>2$ \citep{HA12, khorBook18}.
The Tucker format has stable algorithms \citep{Khoromskij_Low_Tacker}, 
but the storage and complexity costs are 
 $\C{O}(d\,r\,n +r^d)$, i.e.\ they grow exponentially with $d$. 
The TT format is a bit more complicated, but does not have this disadvantage. Applications to UQ problems are described in \cite{dolgov2015polynomial, DolgovLitvLiu19, dolgov2014computation}.

The higher the order of the tensor, the more possibilities there are to find
a low-rank approximation \citep{HA12}, as such
tensors   contain not only rows and columns, but also
\emph{slices} and \emph{fibres} \citep{Kolda:01, Kolda:07, DMV-SIAM2:00, HA12}. 
These can be analysed for linear dependencies, super symmetry,
or sparsity,  and may result in a strong data compression \cite{WhyLowRank19}.  


%
%
\subsection{The canonical polyadic (CP) tensor format}   \label{SS:CPFormat}
The CP representation of a multivariate function is one of the easiest tensor representations, 
it was developed in \citep{Hitch:27}. 
%
A tensor representation is a multi-linear mapping. There are many different representations, e.g.\ the
CP tensor representation maps a tensor $\tnb{w}$ into a sum with $r$ terms 
$\sum_{i=1}^r \Ten_{\nu=1}^d \vek{w}_{i,\nu}$. The number $r$ is called the tensor rank.
The storage required to store a tensor in the CP tensor format is  $\C{O}(r\,d\,n)$ 
(assuming that $n=n_1=\ldots=n_d$).  
%
%
%
Denoting the set of all rank-$r$ tensors by $\C{T}^r$, it is evident that if 
$\tnb{w}_1,\tnb{w}_2 \in \C{T}^r$, then the sum  $\tnb{w}_1+\tnb{w}_2 \notin  \C{T}^r $, 
and, therefore, $\C{T}^r$ is not a vector space \citep{Espig2012VarCal}. 
This sum is generally in $\C{T}^{2r}$.


A complete description of fundamental operations in the canonical
tensor format and their numerical cost can be found in
\citep{HA12}. For recent algorithms in the canonical tensor format we
refer to \citep{hackbusch2012tensor,  ModelReductionBook15, EspigDiss, EspigGrasedyckHackbusch2009, ESHA09_1}.
%
To multiply a tensor $\tnb{w}$ by a scalar $\alpha\in\D{R}$ costs $\C{O}(r\, n\, d)$:
\begin{equation}  
\label{eq:CP-scal-prod}
 \alpha\cdot \tnb{w} = \sum_{j=1}^{r} \alpha \bigotimes_{\nu=1}^d \vek{w}_{j,\nu}=
\sum_{j=1}^{r}  \bigotimes_{\nu=1}^d {\alpha_{\nu}}\vek{w}_{j,\nu},
\end{equation}
where $\alpha_{\nu} := \sqrt[d]{|\alpha|}$ for all $\nu > 1$,
and $\alpha_{1} := \sign(\alpha)\sqrt[d]{|\alpha|}$.

The sum of two tensors costs
only $\C{O}(1)$:
\begin{equation}  \label{eq:CP-add}
\tnb{w} = \tnb{u} + \tnb{v} =\left(\sum_{j=1}^{r_u} 
   \bigotimes_{\nu=1}^d \vek{u}_{j,\nu}\right) +
    \left(\sum_{k=1}^{r_v} \bigotimes_{\mu=1}^d \vek{v}_{k,\mu}\right)
  =
  \sum_{j=1}^{r_u+r_v} \bigotimes_{\nu=1}^d \vek{w}_{j,\nu} ,
\end{equation}
where $\vek{w}_{j,\nu}:=\vek{u}_{j,\nu}$ for $j\leq r_u$ and 
$\vek{w}_{j,\nu}:=\vek{v}_{j,\nu}$ for $r_u< j\leq r_u+r_v$. 
The sum $\tnb{w}$ generally has rank $r_u+r_v$.
The Hadamard product can be written as follows
\begin{equation}
\label{eq:CPHadamard}
\tnb{w}=\tnb{u} \odot \tnb{v} =
  \left(\sum_{j=1}^{r_u} \bigotimes_{\nu=1}^d \vek{u}_{j,\nu}\right) 
\odot 
  \left(\sum_{k=1}^{r_v} \bigotimes_{\nu=1}^d \vek{v}_{k,\nu}\right)
=
\sum_{j=1}^{r_u}\sum_{k=1}^{r_v} \bigotimes_{\nu=1}^d \left(\vek{u}_{j,\nu} \odot
\vek{v}_{k,\nu}\right).
\end{equation}
The new rank can increase till $r_ur_v$, and the computational cost is 
$\C{O}(r_u \,r_v n\, d)$.
The scalar product can be computed as follows:
\begin{equation}
\label{eq:CP_scalarProd}
    \bkt{\tnb{u}}{\tnb{v}}_{\C{T}} =\bkt{\sum_{j=1}^{r_u} \bigotimes_{\nu=1}^d \vek{u}_{j,\nu}}%
{\sum_{k=1}^{r_v} \bigotimes_{\nu=1}^d \vek{v}_{k,\nu}}_{\C{T}}
=
\sum_{j=1}^{r_u}\sum_{k=1}^{r_v} \prod_{\nu=1}^d 
\bkt{\vek{u}_{j,\nu}}{\vek{v}_{k,\nu}}_{\C{P}_\nu}.
\end{equation}
The computational cost is $\C{O}(r_u\,r_v\,n\,d)$. 
The tensor rank can be truncated via the ALS-method or Gauss-Newton-method \citep{EspigDiss, Espig2012VarCal}. The scalar product above helps to compute the Frobenius norm $\Vert \tnb{u} \Vert_2:=\sqrt{\bkt{\tnb{u}}{\tnb{v}}_{\C{T}}}$.




\subsection{The Tensor Train format}
\label{app:TT}
The tensor train (TT) format is described in \citep{OSTY09, oseledets2011, HA12,  khorBook18}.
As already noted, it was originally developed in quantum chemistry as ``matrix product
states'' (MPS), see  \citep{VI03} and references therein, and rediscovered later
\cite{oseledetsTyrt2010, oseledets2011}.

\begin{definition}[TT-Format, TT-Representation, TT-Ranks]\label{D:TTFormat}
The \emph{TT-tensor format} is for variable \emph{TT-representation
ranks} $\vek{r}=(r_0, \dots, r_{d})\in \D{N}^{d+1}$ --- with $r_0 = r_d = 1$ and
under the assumption that $d>2$ --- defined by the
following multilinear mapping
\begin{align}   \label{eq:TTRepresentation}
  \tns{U}_{\TT} : \C{P}_{\TT,\vek{r}} &:=  \bigtimes_{\nu=1}^{d} 
  \C{P}_\nu^{r_{\nu-1} \times r_{\nu}}
  \rightarrow   \C{T},\quad \C{P}_\nu = \D{R}^{M_\nu}\; (\nu=1,\dots,d),   \\ \nonumber
   \C{P}_{\TT,\vek{r}}\ni\tnb{P} &=(\tnb{W}^{(\nu)} = (\vek{w}^{(\nu)}_{j_{\nu-1} j_\nu})
   \in\C{P}_\nu^{r_{\nu-1} \times r_{\nu}} :1\le j_\nu\le r_\nu, 1\le\nu\le d) \\ \nonumber
   &\qquad\mapsto 
   \tns{U}_{\TT}(\tnb{P}):=\tnb{w}=\sum_{j_0=1}^{r_0}\dots\sum_{j_d=1}^{r_d} 
     \Ten_{\nu=1}^{d} \vek{w}_{j_{\nu-1} j_\nu}^{(\nu)} \; \in \C{T}.
\end{align}
We call  $\tnb{w}:=(\tns{w}_{i_1\dots i_d})=U_{\TT,\vek{r}}(\tnb{P})$ a tensor 
represented in the train tensor format.
Note that the TT-cores $\tnb{W}^{(\nu)}$ may be viewed as a vector 
valued $r_{\nu-1}\times r_\nu$ matrix
with the vector $\vek{w}^{(\nu)}_{j_{\nu-1} j_\nu}\in\D{R}^{M_\nu}$ with the components
$w^{(\nu)}_{j_{\nu-1} j_\nu}[i_\nu]: 1\le i_\nu \le M_\nu$
at index position $(j_{\nu-1}, j_\nu)$.    The representation in components is then
\begin{equation}   \label{eq:TTRep-comp-1}
  (\tns{w}_{i_1\dots i_d})=\sum_{j_0=1}^{r_0} \dots \sum_{j_{d}=1}^{r_{d}} 
  w_{j_0 j_1}^{(1)}[i_1] \cdots 
  w_{j_{\nu-1} j_\nu}^{(\nu)}[i_{\nu}] \cdots w_{j_{d-1} j_d}^{(d)}[i_d]
\end{equation}
Alternatively, each TT-core $\tnb{W}^{(\nu)}$ may be seen as a vector
of $r_{\nu-1}\times r_\nu$ matrices $\vek{W}^{(\nu)}_{i_\nu}$ of length $M_\nu$, i.e.\
$\tnb{W}^{(\nu)} = (\vek{W}^{(\nu)}_{i_\nu}):1\le i_\nu \le M_\nu$.  Then
the representation \refeq{eq:TTRep-comp-1} reads
\begin{equation}   \label{eq:TTRep-comp-2}
  (\tns{w}_{i_1\dots i_d})= \prod_{\nu=1}^d \vek{W}^{(\nu)}_{i_\nu},
\end{equation}
which explains the name \emph{matrix product state}.  Observe that the first matrix
is always a row vector as $r_0=1$, and the last matrix is always a column vector as $r_d=1$.
The matrix components $\vek{W}^{(\nu)}_{i_\nu}$ of the TT-cores are also called ``carriages''
or ``waggons'' with ``wheels'' $i_\nu$ at the bottom, coupled to the next ``carriage''
or ``waggon'' via the matrix product.  This explains the \emph{tensor train} name.
If one notes more carefully $\tnb{W}^{(\nu)} \in \C{P}_\nu^{r_{\nu-1} \times r_{\nu}} = 
\D{R}^{r_{\nu-1}} \otimes \D{R}^{M_\nu} \otimes \D{R}^{r_\nu}$, then \refeq{eq:TTRep-comp-2}
can be written more concisely as
\begin{equation}   \label{eq:TTRep-comp-3}
  \tnb{w} = \tns{U}_{\TT}(\tnb{P}) =  \tnb{W}^{(1)} \times_3^1 \tnb{W}^{(2)} \times_3^1 \cdots
  \times_3^1 \tnb{W}^{(d)} ,
\end{equation}
where $\tnb{U} \times_k^\ell \tnb{V}$ is a contraction of the $k$-th index of $\tnb{U}$
with the $\ell$-th index of $\tnb{V}$, where one often writes just $\times_k$ for
$\times^1_k$.  Thus in \refeq{eq:TTRep-comp-3} the contractions leave the indices from
the $\D{R}^{M_\nu}$ untouched, so that the tensor $\tnb{w}$ is formed.
\end{definition}

Each TT-core (or \emph{block}) $\tnb{W}^{(\nu)}$ is defined by 
$r_{\nu-1}\times r_{\nu} \times M_\nu$ numbers. Assuming $n=M_\nu$ for all $\nu=1,\ldots,d$,
the total number of entries scales as $\mathcal{O}(d\,n\,r^2)$, which is tractable 
as long as $r=\max\{r_k\}$ is moderate.

We follow to the work of Oseledets \citep{oseledets2011} and list the major properties of the 
TT-tensor format.

\paragraph{The multiplication with scalar} $\alpha$ could be simply
done by multiplying one of the TT-cores $\tnb{W}^{(\nu)}$ in the
representation \refeq{eq:TTRep-comp-3} for any $\nu$ in
$\tnb{w} =  \tnb{W}^{(1)} \times_3^1 \tnb{W}^{(2)} \times_3^1 \cdots
  \times_3^1 \tnb{W}^{(d)}$.
But to balance the effect better, define $\alpha_{\nu} := \sqrt[d]{|\alpha|}$ for all $\nu > 1$,
and $\alpha_{1} := \sign(\alpha)\sqrt[d]{|\alpha|}$.  Then $\Ttil{w} = \alpha\cdot\tnb{w}$
is given by
\[  \Ttil{w} =  (\alpha_1\cdot\tnb{W}^{(1)}) \times_3^1 
  (\alpha_2\cdot\tnb{W}^{(2)}) \times_3^1 \cdots  \times_3^1 (\alpha_d\cdot\tnb{W}^{(d)})
  = \Ttil{W}^{(1)}\times_3^1 \cdots \times_3^1 \Ttil{W}^{(d)}.
\]
The new cores are given by $\Ttil{W}^{(\nu)} = (\vtil{W}^{(\nu)}_{i_\nu})
= (\alpha_\nu \vek{W}^{(\nu)}_{i_\nu})$, a sequence of new ``carriage'' matrices. 
The computational complexity is $\mathcal{O}(d\,n\,r^2)$.

\paragraph{Addition of two TT-tensors}
Assume two tensors $\tnb{u}$ and $\tnb{v}$ are given in the TT-tensor format
as in \refeq{eq:TTRep-comp-2}, i.e.\
$(\tns{u}_{i_1\dots i_d})= \prod_{\nu=1}^d \vek{U}^{(\nu)}_{i_\nu}$ and 
$(\tns{v}_{i_1\dots i_d})= \prod_{\nu=1}^d \vek{V}^{(\nu)}_{i_\nu}$.  The sum
$\tnb{w}=\tnb{u}+\tnb{v}$ is given by the new cores $\vek{W}^{(\nu)}_{i_\nu}$ such that
$(\tns{w}_{i_1\dots i_d})= \prod_{\nu=1}^d \vek{W}^{(\nu)}_{i_\nu}$, where
\begin{equation*}
    \vek{W}^{(\nu)}_{i_\nu}=
\begin{pmatrix}
\vek{U}^{(\nu)}_{i_\nu} & \vek{0} \\
 \vek{0} & \vek{V}^{(\nu)}_{i_\nu} 
\end{pmatrix}, \quad 1\le i_\nu \le r_\nu, 2\le\nu\le d-1.
\end{equation*}
and the first and the last cores will be 
\[
\vek{W}^{(1)}_{i_1}=\left(\vek{U}^{(1)}_{i_1}\; \vek{V}^{(1)}_{i_1} \right)  \quad \text{ and } 
\quad
\vek{W}^{(d)}_{i_d}=\left( \begin{array}{c} \vek{U}^{(d)}_{i_d} \\ \vek{V}^{(d)}_{i_d}
\end{array} \right).
\]
As only storage may have to be concatenated, the computational cost is $\C{O}(1)$,
but as  the carriages resp.\ TT-cores grow, the final rank will generally be the
sum of the ranks.

\paragraph{The Hadamard product} $\tnb{w}=\tnb{u} \odot \tnb{v}$ in the TT format
 is computed as follows.  Assume two tensors $\tnb{u}$ and $\tnb{v}$ are given in 
 the TT tensor format as in \refeq{eq:TTRep-comp-2}, i.e.\
$(\tns{u}_{i_1\dots i_d})= \prod_{\nu=1}^d \vek{U}^{(\nu)}_{i_\nu}$ and 
$(\tns{v}_{i_1\dots i_d})= \prod_{\nu=1}^d \vek{V}^{(\nu)}_{i_\nu}$.  The Hadamard product is
\[
    (\tns{w}_{i_1\dots i_d}) = (\tns{u}_{i_1\dots i_d} \cdot \tns{v}_{i_1\dots i_d}) .
\]
The tensor $\tnb{w}$ has also the TT-tensor format, namely with the new cores
\[
   \vek{W}^{(\nu)}_{i_\nu} = \vek{U}^{(\nu)}_{i_\nu}\otimes_{\mrm{K}} \vek{V}^{(\nu)}_{i_\nu}, 
   \quad 1\le i_\nu \le r_\nu, 1\le\nu\le d ,
\]
where $\otimes_{\mrm{K}}$ is the Kronecker product of two matrices \citep{HA12}.
The rank of $\tnb{W}^{(\nu)} = (\vek{W}^{(\nu)}_{i_\nu})$ is the product
of the ranks of the TT-cores  $\tnb{U}^{(\nu)}$ and $ \tnb{V}^{(\nu)}$.

\paragraph{The Euclidean inner product} of two tensors in the TT-format
as in \refeq{eq:TTRepresentation} 
\begin{equation*}
\tnb{u} = \sum_{j_0=1}^{r^u_0}\dots\sum_{j_d=1}^{r^u_d} 
     \Ten_{\nu=1}^{d} \vek{u}_{j_{\nu-1} j_\nu}^{(\nu)} , \qquad
\tnb{v} = \sum_{j_0=1}^{r^v_0}\dots\sum_{j_d=1}^{r^v_d} 
     \Ten_{\nu=1}^{d} \vek{v}_{j_{\nu-1} j_\nu}^{(\nu)} ,
\end{equation*}
with ranks $\vek{r}^u$ and $\vek{r}^v$ can be computed as follows:
\[
\bkt{\tnb{u}}{\tnb{v}}_{\C{T}}= \sum_{j_0=1}^{r^u_0}\dots\sum_{j_d=1}^{r^u_d}
  \sum_{i_0=1}^{r^v_0}\dots\sum_{i_d=1}^{r^v_d}  \prod_{\nu=1}^d \;
       \bkt{\vek{u}_{j_{\nu-1} j_\nu}^{(\nu)}}{\vek{v}_{i_{\nu-1} i_\nu}^{(\nu)}}_{\C{P}_\nu}.
\]
The computational complexity is $\C{O}(d\,n\,r^4)$, and can be reduced further \citep{OselTT}. 
%
%
%

\paragraph{Rank truncation in the TT format.}
The rank truncation operation is based on the SVD algorithm and requires 
$\C{O}(d\,n\,r^3)$ operations \citep{GrasHack11}.  The TT-rounding algorithm (p. 2305 in \citep{oseledets2011}) is based on QR decomposition and costs $\C{O}(d\,n\,r^3)$.

Corollary 2.4 in \citep{oseledets2011} states that for a given tensor $\tnb{w}$ and
rank bounds $r_k$, the best approximation to $\tnb{w}$ in the Frobenius norm with TT-ranks 
bounded by $r_k$ always exist (denote it by $\tnb{w}^*$), and the TT-approximation 
$\tnb{u}$ computed by the TT-SVD algorithm (p. 2301 in \citep{oseledets2011}) is quasi-optimal:
\begin{equation}
    \Vert \tnb{w}-\tnb{u} \Vert_F\leq \sqrt{d-1} \Vert \tnb{w} - \tnb{w}^{*} \Vert_F.
\end{equation}
In \cite{Kressner17_Trunc} the authors suggested a new re-compression 
randomised algorithm for Tucker and TT tensor formats. 
For the rank-adaptive DMRG-cross algorithm see \cite{so-dmrgi-2011proc}, and its extension in \cite{Dolgov14_AmenCross}.

\section{Algorithms} \label{S:algs}
%
As one may have gleaned from the preceding,
very often not only the density $p_{\vX}(\by)$ is of interest, but 
expected values of functions of the density, 
e.g.\ in order to compute an $f$-divergence or the entropy. 
{This section lists well-known numerical algorithms to actually compute functions 
from Table~\ref{table:f-divergence}} like
%
$f(p_{\vX}(\by))$ by $f(\tnb{P})$, where the $f$'s considered are
\begin{equation}   \label{eq:nec-fncts}
f(\cdot)=\{\sign(\cdot), (\cdot)^{-1}, \sqrt{\cdot}, \sqrt[m]{\cdot}, 
  (\cdot)^k, \log{(\cdot)},\exp{(\cdot)},(\cdot)^2, \vert \cdot \vert\},
\end{equation}
$k>0$, and
$\tnb{P}$ is a tensor which represents the values of $p_{\vX}(\by)$ on a
discretisation grid as explained in \refSS{low-rank-RV-th}, i.e.\
$\tnb{P} = p_{\vX}(\That{X})= \sum_{j=1}^{r_p} \Ten_{\nu=1}^d \vek{p}_{j,\nu}$.
The $m$-th root is needed for scaling.

Here we collect iterative algorithms and series expansions which have appeared in various places 
scattered in the literature, and which are used to compute functions in matrix algebras.
As was shown in \refS{theory}, these can also be used in the Hadamard algebra,
which is isomorphic to an algebra of diagonal matrices.
These algorithms are non-trivial and their detailed discussion is beyond the scope of
this work, for this we point to the relevant literature. Specific possible difficulties which may 
appear if they are used in the low-rank representaion are:
the intermediate tensor ranks of iterates may become very large, or a stable rank truncation 
procedure may either not be available or it may be very computationally intensive. 
Included are also algorithms specific to low-rank tensor representations,
the so-called \emph{cross-approximation}, see \refSS{TT-cross-apr}.


\subsection{Consistency of data}  \label{SS:consistency}
Making sure the data one works with is consistent is of great importance, and
should be checked before any other computations.  This was already touched upon at the
end of \refSS{low-rank-RV-th}.   Recall that for the continuous case, the \pdf has to 
be in the set $\F{D}$ (cf.\ \refeq{eq:set-prob-dens-fct}), and the \pcf has to be in 
the set $\F{C}$ (cf.\ \refeq{eq:set-char-fct}).  The discussion at the
end of \refSS{low-rank-RV-th} also showed that these conditions are most easily checked
by making sure that the discretised \pdf  (cf.\ \refeq{eq:pdf-tensor-1}) 
is real valued and positive $\tnb{P} \ge \tnb{0}$, and that the discrete integral 
\refeq{eq:discr-int} evaluates to unity, $\C{S}(\tnb{P}) = 1$.

That the \pdf is real valued is translated by the Fourier transform (FT) 
into the \pcf being Hermitean.  The discrete FT preserves these
properties, and  these properties are preserved also
in the low-rank representation.  
To check that the discretised low-rank representation of the \pdf is real valued,
one only has to make sure that the imaginary part vanishes---it does not even
have to be represented.  

The condition that the \pdf is non-negative, which the Fourier transform (FT)
translates into the \pcf being positive definite, is preserved also by the discrete FT.
But this is a local, point-wise condition, and it may be violated in the process
of compression or low-rank approximation, and is also not so easily checked in
a compressed low-rank format.
It was already stated in \refSS{low-rank-RV-th} that this is equivalent to the
positivity of the convolution operator $\tnb{K}_{\tnb{\Phi}}$ using the \pcf tensor 
$\tnb{\Phi}$, whose spectrum are the values of the \pdf tensor $\tnb{P}$,
but this is not easily checked either.  As will be shown later in \refSS{itrunc}, 
it is possible to compute the level set function $\Lambda_U(\tnb{P})$ \refeq{eq:level-set-fct} 
for the level set $U=]-\infty, 0]$.  In case this happens to be non-zero, indicating some indices 
$\vek{i}$ where $\tns{P}_{\vek{i}} \ge 0$ may be violated, one obtains a corrected 
\pdf tensor $\tnb{P}_c$ by seting 
\begin{equation}  \label{eq:corr-P-tens-neg}
   \tnb{P}_c := \tnb{P} - \Lambda_U(\tnb{P}).
\end{equation}
This new approximation has all the faulty negative values removed, so that
$\tnb{P}_c \ge \tnb{0}$.  This change will not affect the Hermitean or self-adjoint
character discussed earlier.

Through a correction like \refeq{eq:corr-P-tens-neg}, or just due to the the process
of compression or low-rank approximation, the condition $\C{S}(\tnb{P}) = 1 =\beta$ may not be
satisfied any more.  In case this is not satisfied, say $1 \ne\beta >0$, 
for a re-scaled $\tnb{P}_{s}$ one may set
\begin{equation}  \label{eq:corr-P-tens-one}
   \tnb{P}_{s} := \tnb{P} / \beta.
\end{equation}
This re-scaling will not affect the positive character of $\tnb{P}$.

\subsection{Overview of methods to compute tensor functions}  \label{SS:comp-tens}
The topic of point-wise functions was already discussed at the
end of \refSS{low-rank-RV-th}.
There are many methods to compute $f(\tnb{w})$. Most of these only use the fact that the 
underlying structure (here of $\C{T}$) is that of a unital  C*-algebra \citep{Segal1978},
and---at least implicitly---define the functions via spectral calculus.
Below we have collected and adapted some of the classical matrix algorithms
for computing the tensor functions listed above in \refeq{eq:nec-fncts}.


Many algorithms developed (cf.\ \citep{NHigham}) to compute functions of matrices use
only operations from the underlying algebra.  They can thus be used in any unital C*-algebra.
Exponential sums to
compute ${\tnb{w}}^{\odot -1}$ 
are implemented
in \citep{BraessHack05, Beylkin05}.  Some other methods and estimates are developed 
in \citep{Grasedyck2004, Hack06_Kron_nonlocal, Khor06_StruckRank}.
The function ${\tnb{w}}^{\odot -1/2}$ was
computed in \citep{khorBook18, Khor06_StruckRank}, and
 ${\tnb{w}}^{\odot -\mu}$ in \citep{KhorTensors07, Khor06_StruckRank}.
A quadrature rule to compute the Dunford-Cauchy contour integral for holomorphic fucntions 
mentioned at the end of \refSS{low-rank-RV-th}
is presented in \citep{GAHAKH03, GAHAKH05, HKTP, LowKroneker, litv17Tensor, khorBook18}.
Iterative methods of Newton and Newton-Schultz type are used in
\citep{EspigDiss, Espig2012VarCal, TensorCalEspig, Espig:2013, ESHALIMA11_1}. 
The inversion ${\tnb{w}}^{\odot -1}$ in the Tucker format is done by the
Newton–Schultz approximate iteration in \citep{Oseledets2009}.

The various TT-cross approximation algorithms \citep{oseledetsTyrt2010, BallaniGrasedyck2015, 
ds-alscross-2017pre, DolgovLitvLiu19}
are a well-known alternative to iterative methods and series expansions, and they are sketched later 
in \refSS{TT-cross-apr}. 
Whereas the previously mentioned algorithms work in any unital C$^*$-algebra, the cross approximation
algorithms make use of a specific representation.  We focus here on the tensor train (TT) representation.
The TT-cross approximation computes a low-rank TT approximation of $f(\tnb{w})$ for a given 
function $f(\cdot)$ and a tensor $\tnb{w}$ directly, ``on the fly''.
The TT-cross algorithms assume that the full tensor $\tnb{v}:=f(\tnb{w})$ is not given explicitly, 
but rather as a function which can return any element $\tnb{v}_{\tnb{i}}:=f(\tnb{w}_{\tnb{i}})$ 
for a given index $\tnb{i}$.
The Tucker cross algorithm to approximate $\tnb{w}^{\odot -1}$
was used in \citep{OselSavoTyrt08}.  In \citep{Savo08}, the multigrid Cross 3D
algorithm is used to calculate the Gauss polynomial. In \citep{Osel10_Cross},
iterative methods in the Tucker format are used to compute $ \tnb{w}^{\odot 1/3}$. 
In \citep{BallaniGrasedyck2015, EspigGrasedyckHackbusch2009}, the authors are 
using the cross method for the Hierarchical
Tucker format to estimate various functionals of the solution in the UQ context.
In \citep{HeidelKhor18, litv17Tensor} the sinc quadrature is used to compute fractional
derivatives of the $d$-dimensional Laplace operator.
The TT-approximation of the sign function was computed in
\citep{dolgov2015polynomial, dolgov2014computation}, but the TT-ranks were large.
The point-wise inverse, level sets, and sign functions are defined and computed in \citep{ESHALIMA11_1}.
The point-wise inverse is required for computing the point-wise functions
$\sign(\tnb{w})$ and $\sqrt{\tnb{w}}$, \citep{ESHALIMA11_1}. The Newton algorithm for 
computing $\sign(\tnb{w})$ and $\exp(-\tnb{w})$, its convergence and error analysis for 
hierarchical matrices were considered in \citep{Grasedyck03_Riccati}.

%
%
%

 %
Polynomials are certainly the simplest functions, as they are directly elements of the algebra.  
General polynomials one would numerically
compute by e.g.\  using Horner's scheme in order to minimise the number
of Hadamard products.  Here only simple powers will be needed
as an auxiliary function 
for a  kind of ``scaling'', 
namely iterated squaring or iterated inverse squaring.  
In some of the algorithms the scaling factor involves e.g.\ $\nd{\tnb{w}}_\infty$, 
we shall sketch at the end of the following \refSS{itrunc} how to compute that. 

Denoting the power function by $\tnb{w} \mapsto \tnb{w}^{\odot m}=\tns{\Psi}_{\text{pow}}(m,\tnb{w})$,
one also wants to use it for negative powers; for  $m < 0$ this is simply 
$\tns{\Psi}_{\text{pow}}(m,\tnb{w}) =\tns{\Psi}_{\text{pow}}(-m,\tnb{w}^{\odot -1})$,
and $\tnb{w}^{\odot -1}$ is shown in \refSS{itrunc}.
A really simple and well known way to 
make sure that no unnecessary multiplications are performed is for $m \in\D{N}$ given by the recursive
formula (although usually implemented in a loop):
\begin{equation}  \label{eq:pow-m}
  \tns{\Psi}_{\text{pow}}(m,\tnb{w}) = 
  \begin{cases}
     m>1 \text{ and odd}:  & \tnb{w} \odot \tns{\Psi}_{\text{pow}}(m-1,\tnb{w}); \\
     m  \text{ even}: &  \tns{\Psi}_{\text{pow}}(\frk{m}{2},\tnb{w}) \odot 
                                \tns{\Psi}_{\text{pow}}(\frk{m}{2},\tnb{w}); \\
     m = 1: &  \tnb{w}; 
  \end{cases}
\end{equation}


\ignore{
\begin{algorithm}
  \caption{Computing $\tns{\Psi}_{\text{pow}}(m,\tnb{w})$}
       \label{alg:pow-m}
  \begin{algorithmic}[1]
     \Require $(m \ge 0)\; \vee \; (\tnb{w}^{\odot -1} \text{exists}$).
     \Ensure Output $\tnb{v} = \tnb{w}^{\odot m}$
     \State $\tnb{v} \gets \tnb{1}$
     \If {$m < 0$}
       \State $\tnb{x} \gets T_{\epsilon}(\tnb{w}^{\odot -1}$);
       \State $n \gets -m $
     \Else
       \State $\tnb{x} \gets \tnb{w}$
       \State $n \gets m$
     \EndIf
     \While {$n > 0$} 
       \If{$n$ is even}
          \State $\tnb{x} \gets T_{\epsilon}(\tnb{x} \odot \tnb{x})$
          \State $n \gets n/2$
       \Else   
          \State $\tnb{v} \gets T_{\epsilon}(\tnb{x} \odot \tnb{v})$
          \State $n \gets n-1$
       \EndIf
     \EndWhile
  \end{algorithmic}
\end{algorithm}
}

\begin{remark}
\label{rem:Rank-trunc}
While performing algebraic operations on tensors in some compressed format, the tensor ranks are increasing. To keep computational cost low, a tensor compression 
after each or after a number of algebraic operations
may be necessary. 
We perform the \emph{truncation}
$\tns{T}_\epsilon$ to \emph{low rank} $r$ with error $\epsilon$
\citep{HA12, Espig:2013}.
We thus allow that the algebraic operations are possibly
only executed approximately.  It is assumed that such an approximation is performed whenever necessary, 
and it will not be always pointed out explicitly.
\end{remark}

%

%

%
%

%
\subsection{Iterative methods} \label{SS:itrunc}
We want to compute $f(\tnb{w})$ for some function $f:\C{T}\to\C{T}$ from the list at the beginning of this section. We describe how to do it through iteration (see also \cite{Matthies2020}).  Thus we have an iteration function $\tns{\Psi}_f$, 
which only uses operations from the Hadamard algebra on $\C{T}$, and which is iterated, 
$\tnb{v}_{i+1}=\tns{\Psi}_f(\tnb{v}_i)$, and converges to a fixed point $\tns{\Psi}_f(\tnb{v}_*)=\tnb{v}_*$. 
When started with a $\tnb{v}_0$ depending on $\tnb{w}$, the fixed point is $\lim_{i\to\infty} \tnb{v}_i = \tnb{v}_* =\tns{\Psi}_f(\tnb{v}_*)= f(\tnb{w})$.
A Newton-type family of high-order iterative methods for some matrix functions
was discussed in \cite{Newton18}.


%
To deal with truncation in iterative algorithms,
the standard iteration map
$\tns{\Psi}_f$ is replaced by $\tns{T}_\epsilon \circ \tns{\Psi}_f$.
Here one speaks about a \emph{perturbed or truncated iteration} \citep{whBKhEET:2008,
Zander10, Matthies2020}.
The general structure of the iterative
algorithms for a post-processing task $f$ 
or for an auxiliary function is shown in Algorithm~\ref{alg:basic}.
\begin{algorithm}
  \caption{Iteration with truncation}
       \label{alg:basic}
   \begin{algorithmic}[1] 
    \State Input: tensor $\tnb{v}_0$; output: tensor $\tnb{v}_n$ after $n$ iterations
    \State Start with some initial \emph{compressed} guess  $\tnb{v}_0=\tnb{w}$.
    \State $i\gets 0$
    \While  {\emph{no convergence}} 
       \State $\tnb{v}_{i+1} \gets \tns{T}_\epsilon \circ \tns{\Psi}_f(\tnb{v}_i)$
       \State $i\gets i+1$
    \EndWhile 
  \end{algorithmic}
\end{algorithm}
In case the iteration by $\tns{\Psi}_f$ is \emph{super-linearly}
   convergent, the \emph{truncated iteration} $\tns{T}_\epsilon \circ \tns{\Psi}_f$ will
   still \emph{converge} super-linearly, but finally \emph{stagnate} in an 
   $\epsilon$-neighbourhood of the fixed point $\tnb{v}_*$ \citep{whBKhEET:2008}.
   If the iteration by $\tns{\Psi}_f$ is \emph{linearly} convergent with \emph{contraction}
   factor $q$, the \emph{truncated iteration}  $\tns{T}_\epsilon \circ \tns{\Psi}_f$ will
still \emph{converge} linearly, but finally
\emph{stagnate} in an $\epsilon/(1-q)$-neighbourhood
of $\tnb{v}_*$ \citep{Zander10}. Ideas how to choose the starting value $\tnb{v}_{0}$ are given in \cite{Matthies2020}.
%

%
%
\ignore{
\paragraph{Computing $\tnb{w}^{\odot 2^k}$ or $\tnb{w}^{\odot -2^k}$ for $k\in\D{N}$}
{The basic function is really simple, but it will be reused several times:
\begin{equation}  \label{eq:square}
  \tns{\Psi}_2(\tnb{v}) = \tnb{v} \odot \tnb{v} .
\end{equation}
}
{
Then, to compute $\tns{\Psi}_{2^k}$, one executes the Algorithm~\ref{alg:2-exp-k-times}:}
\begin{algorithm}
  \caption{Computing $\tns{\Psi}_{2^k}$}
       \label{alg:2-exp-k-times}
  \begin{algorithmic}[1]
     \State Start with some input point $\tnb{v}_0$ and a $k\in\D{Z}$.
     \Require $k \ge 0$
     \Ensure Output $\tnb{v} = \tnb{v}_0^{\odot 2^k}$
     \If {$k = 0$}
       \State $\tnb{v} \gets \tnb{w}$;
       \Statex\Comment{ Only for completeness, no one should use this part. } 
     \ElsIf {$k>0$}
       \State $\tnb{v} \gets \tnb{w}^{\odot 2}$;
     \Else
       \State $\tnb{v} \gets \tnb{w}^{\odot -1}$;
     \EndIf
     \For{$i\gets 2$ to $k$} 
       \State $\tnb{v} \gets \tns{\Psi}_2(\tnb{v})$;
     \EndFor 
  \end{algorithmic}
\end{algorithm}
{This is not really an iterative algorithm, the ``iteration'' is executed exactly $k$ times. It is needed
for scaling in the computation of $\log$ and $\exp$.  The starting point is $\tnb{v}_0 = \tnb{w}$ for 
the computation of $\tnb{w}^{\odot 2^k}$ for $k\ge 0$, and
in case $\tnb{w}^{\odot -2^k}$ is wanted, one takes $\tnb{v}_0 = \tnb{w}^{\odot -1}$.
}
}
\paragraph{Computing the point-wise inverse $\tnb{w}^{\odot -1}$.}
Apply Newton's method to $\tns{F}(\tnb{x}):=\tnb{w}-\tnb{x}^{\odot -1}$ to
approximate the inverse of a given tensor
$\tnb{w}$, and one obtains \citep{Oseledets2009} the 
{following iteration function $\tns{\Psi}_{\odot-1}$ ---to be used with Algorithm~\ref{alg:basic} 
with the initial iterate $\tnb{v}_0=\alpha \cdot\tnb{w}$ 
to bring $\tnb{v}_0$ close to $\tnb{v}_a=\tnb{1}$}:
\begin{equation}  \label{eq:Newton-inverse}
    \tns{\Psi}_{\odot-1}(\tnb{v})= \tnb{v}\odot (2\cdot\tnb{1} − \tnb{w} \odot \tnb{v}) .
\end{equation}

The iteration converges if the initial iterate $\tnb{v}_0$ 
satisfies $\nd{\tnb{1}−\tnb{w}\odot\tnb{v}_0}_\infty <1$.
A possible candidate for the starting value is $\tnb{v}_0=\alpha \tnb{w}$ with
$\alpha < (1 /\nd{\tnb{w}}_\infty)^2$. 
For such a $\tnb{v}_0$, the convergence initial
condition  $\nd{\tnb{1} - \alpha \tnb{w}^{\odot 2}}_\infty <1$ is always satisfied.
{As the initial iterate was scaled, the fixed point of the iteration is $\tnb{v}_* = (1/\alpha) \cdot \tnb{w}^{\odot -1}$, and
thus the final result is $\tnb{w}^{\odot -1} = \alpha \cdot \tnb{v}_*$.}

{
We did not assume that $\tnb{w}$ is invertible, as the iteration actually computes the pseudo-inverse
(i.e.\ only the non-zero entries are inverted).  It is easily seen from \refeq{eq:Newton-inverse} that
with $\tnb{v}_0=\alpha \tnb{w}$ zero entries in $\tnb{w}$ stay zero during the iteration.
}
\paragraph{Computing point-wise $\sqrt{\tnb{w}}$ resp.\ $(\tnb{w})^{\odot \frk{1}{2}}$.}
{
The Newton iteration for $\tns{F}(\tnb{x}):=\tnb{x}^{\odot  2} - \tnb{w}$
uses the iteration function \refeq{eq:Newton-sqrt-1} 
together with Algorithm~\ref{alg:basic}.
}

\begin{equation} \label{eq:Newton-sqrt-1}
  \tns{\Psi}_{\sqrt{}}(\tnb{v}) = \frac{1}{2} \cdot (\tnb{v} + \tnb{v}^{\odot -1}\odot\tnb{w}) .
\end{equation}
The starting value can be $\tnb{v}_0 = (\tnb{w}+\tnb{1})/2$, other stating values
obtained through scaling are described later.

Unfortunately, \refeq{eq:Newton-sqrt-1} contains an inverse power $\tnb{v}^{\odot -1}$ and can thus
not be computed directly only with operations from the algebra.  One could use the just described
algorithm for the inverse, but the nested iteration is usually not so advantageous.
An alternative is the well known stable inversion free Newton-Schulz iteration \citep{NHigham}, 
which simultaneously computes $\tnb{v}_*^+=\tnb{w}^{\odot \frk{1}{2}}$ 
and its inverse $\tnb{v}_*^- =\tnb{w}^{\odot -\frk{1}{2}}$.
Setting $\tnb{V}_0 = [\tnb{y}_0, \tnb{z}_0] = [\alpha \cdot\tnb{w}, \tnb{1}] \in \C{T}^2$, 
the iteration function is 
best written using the auxiliary function $\tns{A}(\tnb{y},\tnb{z}) = 3\cdot\tnb{1} - \tnb{z}\odot\tnb{y}$:
\begin{equation}  \label{eq:Newton-sqrt}
    \tns{\Psi}_{\sqrt{}}\left(\begin{bmatrix}
       \tnb{y} \\ \tnb{z}
    \end{bmatrix}\right) = \frac{1}{2}\begin{bmatrix}
     \tnb{y}\odot\tns{A}(\tnb{y},\tnb{z}) \\  \tns{A}(\tnb{y},\tnb{z})\odot\tnb{z}
    \end{bmatrix} .
\end{equation}

The iteration converges to $\tnb{V}_* = [\tnb{v}_*^+, \tnb{v}_*^-] = [ \sqrt{\tnb{y}_0}, (\sqrt{\tnb{y}_0})^{\odot -1}] $
if $\nd{\tnb{1}-\tnb{y}_0}_\infty < 1$, which can be achieved with a scaling factor $\alpha < 1/\nd{\tnb{w}}_\infty$.  
As the initial iterate was scaled, the fixed point of the iteration is $\tnb{v}_*^+ = \sqrt{\alpha} \cdot \sqrt{\tnb{w}}$ and 
$\tnb{v}_*^- = (1/\sqrt{\alpha}) \cdot (\sqrt{\tnb{w}})^{\odot -1}$.  Thus the final result is
$\sqrt{\tnb{w}} = (1/\sqrt{\alpha}) \cdot \tnb{v}_*^+$ and
$(\sqrt{\tnb{w}})^{\odot -1} = \sqrt{\alpha} \cdot \tnb{v}_*^-$.

\paragraph{Computing $\tnb{w}^{\odot \frk{1}{m}}$ when $m\in\D{N}$.}
{This function may be needed for logarithmic scaling purposes, see \refSS{series}.}
For the theory and required definitions see Section 7 in \citep{NHigham}.
A new family of high-order iterative methods for the matrix $m$-th root was suggested in \citep{pthroot15}. 
To compute the principal $m$-th root of $\tnb{w}$,
{where it is assumed that $\tnb{w}\ge\tnb{0}$, one considers}
Newton's method for $\tns{F}(\tnb{x}) = \tnb{x}^{\odot m}-\tnb{w}=\tnb{0}$.
The iteration function with $\tnb{v}_0=\tnb{w}$
{corresponding to \refeq{eq:Newton-sqrt-1}} looks like
{
\begin{equation}  \label{eq:Newton-p-root}
    \tns{\Psi}_{m-\text{root}}(\tnb{v}) = \frac{1}{m} \left( (m-1)\cdot\tnb{v} + \tns{\Psi}_{\text{pow}}(1-m,\tnb{v})\odot\tnb{v}_0 \right) .
\end{equation}
}
{If $m\ge2$, this involves a negative power $\tnb{v}^{\odot(1-m)}=\tns{\Psi}_{\text{pow}}(1-m,\tnb{v})$.}
The convergence analysis is rather complicated, see more in Section 7.3, \citep{NHigham}{, but the algorithm converges for all $\tnb{w}\ge\tnb{0}$}.

{
Just as there is a more stable version \refeq{eq:Newton-sqrt} for $m=2$ avoiding inverses, so one has a ``double iteration'' \citep{NHigham} here as well.
It is best written using the auxiliary function $\tns{A}(\tnb{y},\tnb{z}) = (1/m)\cdot ((m+1)\cdot\tnb{1} - \tnb{z})$:
\begin{equation}  \label{eq:Newton-sqrt-p}
    \tns{\Psi}_{m-\text{root}}=\left(\begin{bmatrix}       \tnb{y} \\ \tnb{z}
    \end{bmatrix}\right) = \begin{bmatrix}
     \tnb{y}\odot\tns{A}(\tnb{y},\tnb{z}) \\  \tns{\Psi}_{\text{pow}}(m,\tns{A}(\tnb{y},\tnb{z}))\odot\tnb{z}
    \end{bmatrix} ,
\end{equation}
where $\tnb{y}_i \to \tnb{w}^{\odot -\frk{1}{m}}$ and $\tnb{z}_i \to \tnb{w}^{\odot \frk{1}{m}}$.  
The starting values are $\tnb{V}_0 = [\tnb{y}_0, \tnb{z}_0] = 
[\alpha \cdot\tnb{1}, (\alpha)^m \tnb{w}] \in \C{T}^2$, with 
$\alpha < (\nd{\tnb{w}}_\infty/\sqrt{2})^{-\frk{1}{m}}$.
For scaling purposes it is best used with $m=2^k$.
}

Another way of computing the $m$-th root is Tsai's algorithm \citep{TsaiShiehYates1988, LorinTian2021},
which uses the auxiliary function
$\tns{B}(\tnb{y}) = (2\cdot\tnb{1} + (m-2)\cdot\tnb{y})\odot (\tnb{1} + (m-1)\cdot\tnb{y})^{\odot -1}$:
\begin{equation}  \label{eq:Newton-sqrt-m}
    \tns{\Psi}_{Tsai}=\left(\begin{bmatrix}       \tnb{y} \\ \tnb{z}
    \end{bmatrix}\right) = \begin{bmatrix}
     \tnb{y}\odot\tns{\Psi}_{\text{pow}}(m,\tns{B}(\tnb{y})) \\  \tnb{z}\odot(\tns{B}(\tnb{y}))
    \end{bmatrix} ,
\end{equation}
with starting value $\tnb{V}_0 = [\tnb{w}, \tnb{1}]$.  Then $\tnb{z}_i \to \tnb{w}^{\odot \frk{1}{m}}$.

\paragraph{Computing $\sign(\tnb{w})$.}
The tensor $\sign(\tnb{w}) \in \Tp$ is defined point-wise for all $\vek{i}\in \Ip$ by
\begin{equation} \label{eq:sign}
(\sign(\tnb{w}))_{\vek{i}}:=
                       \begin{cases}
               \phantom{-}1, & \tnb{w}_{\vek{i}} > 0 ; \\
                         -1, & \tnb{w}_{\vek{i}} < 0 ; \\
               \phantom{-}0, & \tnb{w}_{\vek{i}} = 0 .
                       \end{cases}
\end{equation}
The equation to be used for the Newton iteration is
$\tns{F}(\tnb{x}):=\tnb{x}\odot \tnb{x}-\tnb{1}$.  With starting value $\tnb{v}_0=\tnb{w}$
one obtains the following iteration function:
{
\begin{equation}  \label{eq:Newton-sign}
 \tns{\Psi}_{\sign}(\tnb{v})=\frac{1}{2}(\mu\cdot\tnb{v}+\frac{1}{\mu}\cdot\tnb{v}^{\odot-1}),
\end{equation}
}
with $\mu = \nd{\tnb{v}^{\odot-1}}_\infty / \nd{\tnb{v}}_\infty$ (cf.\ Section 8.6 in \citep{NHigham}).  
This method converges to $\sign(\tnb{w})$.

Alternatively, one can rewrite this iteration function without computing $\tnb{v}^{\odot -1}$, namely
\begin{equation} \label{eq:Newton-Schulz-sign}
    \tns{\Psi}_{NS}(\tnb{v})=\frac{1}{2}\cdot\tnb{v} \odot ( 3\cdot \tnb{1}-\tnb{v}\odot \tnb{v} ) ,
\end{equation}
{with starting value $\tnb{v}_0 = \alpha\cdot \tnb{w}$,
where $\alpha = \nd{\tnb{w}^{\odot -1}}_\infty / \nd{\tnb{w}}_\infty$.}
The last formula is called the Newton-Schulz iteration, it has quadratic (local) convergence to 
$\sign(\tnb{w})$ \citep{EspigDiss,ESHALIMA11_1,ESHAHARS11_2,NHigham}.

\paragraph{Computing the absolute value $\ns{\tnb{w}}$}is simple if $\sign(\tnb{w})$ 
\refeq{eq:sign} is available (see above):
Having the $\sign(\cdot)$ function, one can compute the absolute value, 
characteristic function of a set, and the level set function:
\begin{equation}  \label{eq:abs-val-fct}
\ns{\tnb{w}} = \tnb{w}\odot \sign{(\tnb{w})}.
\end{equation}

\paragraph{Computing the characteristic function of a sub-set} of values $I \subset \R$
of a tensor $\tnb{w}$ can be done by using the shifted $\sign$-function \refeq{eq:sign} \cite{ESHA09_1,ESHAHARS11_2,Matthies2020}.
The \emph{characteristic} function of $\tnb{w} \in \Tp$ in an interval $I \subset \R$ is a tensor $\chi_I(\tnb{w}) \in \Tp$.
It is defined for every multi-index $\vek{i}$ point-wise as
\begin{equation}   \label{equ:char}
(\chi_I(\tnb{w}))_{\vek{i}}:=
                        \begin{cases}
                          1, & \tnb{w}_{\vek{i}} \in I ; \\
                          0, & \tnb{w}_{\vek{i}} \nin I .
                        \end{cases}
\end{equation}
%
%
%
%
Let $a, b \in \R$. If $I=(-\infty, b)$, then  
$\chi_{I}(\tnb{w})=\frac{1}{2}(\tnb{1} + \sign(b\tnb{1}-\tnb{w}))$. If $I=(a,+\infty)$, 
then $\chi_{I}(\tnb{w})=\frac{1}{2}(\tnb{1} - \sign(a\tnb{1}-\tnb{w}))$. 
And if $I=[a,b]$, then 
\begin{equation}  \label{eq:char-fct}
\chi_{I}(\tnb{w})=\frac{1}{2} (\sign(b\tnb{1}-\tnb{w}) -
\sign(a\tnb{1}-\tnb{w})).
\end{equation}
%
%

%
%
%

\paragraph{Computing the level set function of a sub-set} of values $I \subset \R$ of a tensor 
$\tnb{w}$ may be accomplished using the characteristic function \refeq{eq:char-fct}
\cite{ESHA09_1,ESHAHARS11_2,Matthies2020}:
\begin{equation}   \label{eq:level-set-fct}
\Lambda_{I}(\tnb{w}) = \chi_{I}(\tnb{w}) \odot \tnb{w} .
\end{equation}

\paragraph{Computing $\nd{\tnb{w}}_\infty = \max \tnb{w}=\vrho(\tnb{w})$ and
$\min \tnb{w}$.}
We need these values in almost every algorithm above and for checking the consistency. For instance, if 
$\min \tnb{w} \ge \tnb{0}$ then $\tnb{w}\ge \tnb{0}$.  It was already pointed out that 
$\nd{\tnb{w}}_\infty = \max \tnb{w}=\vrho(\tnb{w})$ is the largest---by magnitude---eigenvalue of the associated linear
operator $\tnb{L}_{\tnb{w}}: \tnb{v} \mapsto \tnb{w}\odot \tnb{v}$, and thus itertive eigenvalue algorithms
can be used to compute it.  One of the simplest is the
power iteration.  It can be modified in this special case here to greatly increase its convergence speed,
essentially by repeated squaring.
The modified power iteration algorithm is described in
\citep{Grasedyck-gamm, Matthies2020}. The
convergence is \emph{exponential} with the rate $|\lambda_2/\lambda_1|^{2^i}$, where $\lambda_1$ and $\lambda_2$ are the largest and the second largest eigenvalues of $\tnb{L}_{\tnb{w}}$.



To compute $\min \tnb{w}$ (the smallest eigenvalue) or some intermediate eigenvalues of
$\tnb{L}_{\tnb{w}}$ resp.\ $\tnb{w}$ we suggest to use
the well-known shifting or inverse shifting functions \citep{ESHALIMA11_1, Matthies2020}. 
These and similar techniques are well known from eigenvalue
calculations of large / sparse symmetric matrices  
\citep{golubVanloan, parlett98, saad92, ds-watk-07}.

\subsection{Series expansions} \label{SS:series}
\paragraph{Computing $\log(\tnb{w})$.}
{We assume that $\tnb{w} > 0$.}
Various algorithms, improvement ideas, stability issues and tricks to compute 
$\log(\tnb{w})$  are discussed in \citep{HighamMLog01}.
Improved inverse scaling and squaring algorithms for the matrix logarithm were suggested 
later in \citep{Higham12}. We suggest to follow these works for the case when $\tnb{w}$ is a tensor. 
The stability of the matrix arithmetic-geometric mean iterations for computing the matrix 
logarithm is investigated in \citep{CompLog16}.
For the algorithms to work well $\tnb{w}$ has to be close to the identity $\tnb{1}$, 
which can be achieved by taking roots:
for $\lambda >0$ one has $\log\left(\tnb{w}^\lambda\right) = \lambda \log \tnb{w}$.

From matrix calculus it is known \citep{HighamMLog01,Higham12} that one of the ways 
to compute $\log(\tnb{w})$ for $\tnb{w}$ close to the identity 
is to truncate the Taylor series {(radius of convergence  $\nd{\tnb{x}}_\infty < 1$)}:
\[
\log(\tnb{1}-\tnb{x})= - \sum_{n=1}^\infty \frac{1}{n} \cdot \tnb{x}^{\odot n}
\]
where $\tnb{x}:=\tnb{1}-\tnb{w}$.  If $\tnb{w}$ is not near to the identity, 
then one may use the relation $\log(\tnb{w})=2^k\log(\tnb{w}^{\odot 1/2^k})$, 
where $\tnb{w}^{\odot 1/2^k}\to \tnb{1}$ as $k$ increases \cite{kenney1989condition}.

Another way to compute the logarithm of a positive 
$\tnb{w}$ is Gregory's series  \citep{NHigham}, which converges for all $\tnb{w} > 0$.
Setting $\tnb{z} = (\tnb{1}-\tnb{w})\odot (\tnb{1}+\tnb{w})^{\odot -1}$, one has
\begin{equation}
\label{eq:GregorySeries}
   \log \tnb{w} = -2 \sum_{k=0}^\infty \frac{1}{2k+1} 
   \cdot \tnb{z}^{\odot (2k+1)}.
\end{equation}
Obviously, $\tnb{z}$ involves an inverse, but it has to be computed only once.
We note that the tensor ranks (according to Section~\ref{app:TT}) in 
\refeq{eq:GregorySeries} may increase very fast. 

%
%
%

%
\paragraph{Computing $\exp\tnb{w}$.}
One of the standard algorithms using power series together with scaling and squaring is explained in \cite{NHigham} (Chapter 10):
\begin{equation}
\label{eq:expw}
    \tnb{u}_{r,s} = \left(\sum_{k=0}^r \frac{1}{k! s^k} \tnb{w}^{\odot k}\right)^{\odot s}.
\end{equation}
Here $\lim_{r\to\infty}\tnb{u}_{r,s}=\lim_{s\to\infty} \tnb{u}_{r,s}=\exp\tnb{w}$.  
It is of advantage to use $s$ from the series of powers of $2$, $s=1,2,4,\dots,2^k$, then the $s$-th power can be computed by squaring. For the scaling the best choice is $\alpha > \nd{\tnb{w}}_\infty$.

%
%


%

\subsection{Direct approximation by TT-cross-algorithms}
\label{SS:TT-cross-apr}

In this section, we discuss the so-called cross algorithms \cite{oseledetsTyrt2010,BallaniGrasedyck2015,ds-alscross-2017pre,DolgovLitvLiu19} for $\tnb{v}:=f(\tnb{w})$. 
Some extensions of the TT-Cross and ALS-Cross algorithms (e.g.\ the AMEn algorithm) were suggested by 
S.~Dolgov and co-authors, and can be found in \cite{DolgovALSCross}.

A general cross algorithm
computes the following TT-representation,
see Subsection~\ref{app:TT}, the definition \refeq{eq:TTRepresentation}, and in particular
\refeq{eq:TTRep-comp-1}:
\begin{equation}\label{eq:tt}
  f(\tnb{w}) =
  \tnb{v}(\alpha_1,\ldots,\alpha_M)  = 
  \sum_{s_1=1}^{r_1} \sum_{s_2=1}^{r_2}  \cdots \sum_{s_{M-1}=1}^{r_{M-1}}
      \tns{v}^{(1)}_{s_0,s_1}(\alpha_1) \tns{v}^{(2)}_{s_1,s_2}(\alpha_2) 
      \cdots \tns{v}^{(M)}_{s_{M-1},s_M}(\alpha_M).
\end{equation}
For analytic $f(\cdot)$ the TT-ranks often depend only logarithmically on the 
accuracy \cite{khor-rstruct-2006,uschmajew-approx-rate-2013}.

These algorithms are an alternative to iterations and series expansions, they are
tailored specifically to low-rank approximations in a particular tensor format.
They allow one
to compute a low-rank approximation of $\tnb{v}= f(\tnb{w})$ for a given function $f(\cdot)$ and a tensor $\tnb{w}$ directly, ``on the fly''.
We assume here that the full tensor $\tnb{v}:=f(\tnb{w})$ is not given explicitly, but rather as a function 
which can return any element $\tnb{v}_{\tnb{i}}:= f(\tnb{w}_{\tnb{i}})$ of $\tnb{v}$.

For example, the AMEn algorithm can compute the representation  \refeq{eq:tt} using only $\mathcal{O}(dnr^2)$ entries of $\tnb{v}$ and $\mathcal{O}(dnr^3)$ additional arithmetic operations. 
The pseudocode is listed in \cite{DolgovLitvLiu19}.
It is based on the skeleton decomposition (another name is adaptive cross approximation) of a matrix \cite{TyrtyshACA,ACA,bebe-aca-2011}, 
and the \emph{maxvol} algorithms \cite{gostz-maxvol-2010}. The idea of the \emph{maxvol} algorithms is to find a rank-$r$ matrix approximation $A_r$ of a $n\times m$ matrix $A$, the rank of which is $r$. 
The \emph{maxvol} algorithm suggests to select among all $r\times r$ submatrices the one that has the largest volume (determinant). 
The computational complexity is $\mathcal{O}(r(n+m))$, and the approximation error $\Vert A-A_r \Vert_{\infty}\leq (r+1) \sigma_{r+1}$, where $\sigma_{r+1}$ is the $(r+1)$-th singular value of $A$.

\section{Numerical examples} \label{S:Numerics}
%
{This section shows 
the real applicability of tensor techniques addressed in Sections~\ref{S:theory}--\ref{S:algs}.} 
This includes a validation Example~\ref{ex:validation}, where the KLD is computed with a
well known analytical formula and the AMEn$\_$cross algorithm 
\citep{Dolgov14_AmenCross, so-dmrgi-2011proc} from the TT-toolbox for $d$-dimensional Gaussian \pdfs.
Further, in Example~\ref{ex:validation2} the Hellinger distances
are computed with the well known analytical formula as well as with the AMEn$\_$cross algorithm.
To show the approach on a distribution where the \pdf is not known analytically, the  $d$-variate 
elliptically contoured $\alpha$-stable distributions are chosen and accessed via their \pcfs, and
again KL and Hellinger distances  for different value of $d$, $n$, and the parameter $\alpha$ are
computed (Example~\ref{ex:biVarEllip}).

For the numerical tests below we used the Matlab package \emph{TT-Toolbox} \citep{oseledets2011},
which is well known in the tensor community.  To compute $f(\tnb{w})$, we use the 
\textit{alternating optimization with enrichment (AMEn)} method \citep{Dolgov14_AmenCross}, provided 
in the TT-toolbox library.  This is a block cross algorithm with an error-based enrichment. 
It tries to interpolate the function $f(\tnb{w})$ via the error-enriched maxvol-cross method.
All computations are done on a MacBook Pro computer produced in 2018, equipped with a 6-Core 
Intel Core i7 processor running at 2.2 GHz, and 16 GB RAM.
We started by computing the point-wise inverse, squared root, and exponent of a given discretised 
\pdf represented as a TT tensor.  For this, iterative methods, series expansions, and the AMEn 
algorithm \citep{Dolgov14_AmenCross, so-dmrgi-2011proc} were used in order to make sure
that the AMEn method gives the same results as other methods.

The first Example~\ref{ex:validation} is a validation example, where the analytical formula 
for the KLD is known analytically.  This exact value is compared with the approximate KLD 
(denoted by $\widetilde{D}_{\text{KL}}$) for high values of $d$ and $n$.  
One may observe that they are almost the same.
Additionally, the absolute error ($\text{err}_a:=\vert D_{\text{KL}} - \widetilde{D}_{\text{KL}}\vert$)
as well as the relative error ($\text{err}_r:={\vert D_{\text{KL}} - \widetilde{D}_{\text{KL}}\vert}/
{\vert  D_{\text{KL}}  \vert}$)  and the computing times (last row) are shown.

\begin{example}[Validation example 1 --- KLD]
\label{ex:validation}
Consider two Gaussian distributions  $\C{N}_1:=\C{N}(\vek{\mu}_1,\bC_1)$ and 
$\C{N}_2:=\C{N}(\vek{\mu}_2,\bC_2)$, where $\bC_1:=\sigma_1^2\bI$, 
$\bC_2:=\sigma_2^2\bI$,  $\vek{\mu}_1=(1.1 \ldots, 1.1)$ and $\vek{\mu}_2=(1.4, \ldots, 1.4)\in \D{R}^d$,
$d=\{16,32,64\}$, $\bI$ is the identity matrix, and $\sigma_1=1.5$, $\sigma_2 = 22.1$. 
The well known analytical formula is \citep{pardo2018statistical}
\begin{equation}
\label{eq:dKLD}
       D_{\text{KL}}(\C{N}_1\Vert \C{N}_2)=\frac{1}{2}\left( 
       \tr{(\bC_2^{-1}\bC_1)} + (\vek{\mu}_2 -\vek{\mu}_1)^T\bC_2^{-1} (\vek{\mu}_2-\vek{\mu}_1) - d + \log \left( \frac{\vert \bC_2 \vert}{\vert \bC_1 \vert}\right)
    \right).
\end{equation}
After discretisation of the two \pdfs, one obtains the tensors $\tnb{P}$ and $\tnb{Q}$ 
(cf.\ \refS{intro}).  Then the KLD $\widetilde{D}_{\text{KL}}$ is computed as 
in Table~\ref{table:divergences-discr}. The results are summarised in Table~\ref{table:KLD-comput0}. 
\end{example}

\begin{center}
\begin{table}[t]
 \centering
 \caption{$D_{\text{KL}}$ computed via TT tensors (AMEn algorithm) and the analytical 
 formula \refeq{eq:dKLD} for various values of $d$. TT tolerance = $10^{-6}$,
 the stopping difference between consecutive iterations.}
  \label{table:KLD-comput0} 
 \begin{tabular}{llll}
 \hline
$d$     & 16    & 32     & 64  \\ \hline 
$n$     & 2048  & 2048   & 2048   \\ \hline 
$D_{\text{KL}}$ (exact) & 35.08 & 70.16 & 140.32\\
$\widetilde{D}_{\text{KL}}$ & 35.08 & 70.16 & 140.32 \\ 
$\text{err}_a$ &4.0e-7 & 2.43e-5& 1.4e-5 \\ 
$\text{err}_r$ &1.1e-8 & 3.46e-8& 8.1e-8 \\  
comp.time [s]& 1.0&5.0& 18.7   \\ \hline
\end{tabular}
\end{table}
\end{center}

An important ingredient of the KLD computation is the $\log(\cdot)$ function, 
which can be computed 
via the AMEn method, or the Gregory series  \refeq{eq:GregorySeries}.
We observed that the TT-ranks are increasing very fast in the Gregory series, 
and it is not so transparent how and when to truncate them. 
Therefore, we recommend using the AMEn algorithm (implemented in the TT Toolbox) 
for the approximation $\widetilde{D}_{\text{KL}}$ {by computing on the fly the function 
$f(\tnb{P},\tnb{Q})=\tnb{P}\odot(\log(\tnb{P}) - \log(\tnb{Q}))$}. 
The small absolute  ($\text{err}_a$) and relative ($\text{err}_a$) errors show that the 
AMEn method can be used to compute the KLD. 

The next validation test is with the Hellinger distance: 
\begin{example}[Validation example 2 --- Hellinger distance]
\label{ex:validation2}
For the Gaussian distributions from Example~\ref{ex:validation},
the Hellinger distances computed via the AMEn algorithm---denoted 
by $\widetilde{D}_{H}$ ---and the analytical formula \refeq{eq:Hellinger_an} denoted by $D_{H}$)
are compared in Table~\ref{table:Hellinger-comput0}. 
The squared Hellinger distances for two multi-variate Gaussian distributions can be computed 
analytically as follows (see p.51 and p.45 in \citep{pardo2018statistical}):
\begin{equation}
\label{eq:Hellinger_an}
D_{H}(\C{N}_1, \C{N}_2)^2 =  1 - K_{\frk{1}{2}}(\C{N}_1, \C{N}_2),\quad \text{ where}
\end{equation}
%
\begin{equation}
\label{eq:Hellinger_K}
K_{\frk{1}{2}}(\C{N}_1, \C{N}_2) = \frac{\det(\bC_1)^{\frk{1}{4}}\det(\bC_2)^{\frac{1}{4}}}%
{\det\left(\frac{\bC_1 + \bC_2}{2}\right)^{\frk{1}{2}}} 
\cdot \exp \left( -\frac{1}{8}(\vek{\mu}_1 - \vek{\mu}_2)^{\top}\left(\frac{\bC_1 + \bC_2}{2}\right)^{-1}
(\vek{\mu}_1 - \vek{\mu}_2)\right )
\end{equation}
\end{example}
\begin{center}
\begin{table}[t]
 \centering
  \caption{The Hellinger distance $D_{H}$ computed via TT tensors (AMEn algorithm) and the 
 analytical formula \refeq{eq:Hellinger_an} for various values of $d$. TT tolerance = $10^{-6}$. 
 The Gaussian mean values and covariance matrices are defined in Example~\ref{ex:validation}.}
  \label{table:Hellinger-comput0} 
 \begin{tabular}{llll}
 \hline
$d$     & 16    & 32     & 64  \\ \hline 
$n$     & 2048  & 2048   & 2048   \\ \hline 
$D_{H}$ (exact) & 0.99999 & 0.99999 & 0.99999\\ 
$\widetilde{D}_{H}$ & 0.99992 & 0.99999 & 0.99999 \\ 
$\text{err}_a$ & 3.5e-5& 7.1e-5& 1.4e-4 \\ 
$\text{err}_r$ & 2.5e-5& 5.0e-5& 1.0e-4 \\ 
comp.time [s]& 1.7 & 7.5& 30.5   \\ \hline 
\end{tabular}
\end{table}
\end{center}

The results in Table~\ref{table:Hellinger-comput0} show that the AMEn algorithm is able to 
compute the $D_{H}$ Hellinger distance between two multivariate Gaussian distributions for 
large dimensions $d=\{16,32,64\}$, and for large $n=2048$. The exact and approximate 
values are identical, and the error is small. The absolute and relative errors 
($\text{err}_a$, $\text{err}_r$) can be further decreased by taking a smaller TT tolerance.



After these validation tests, we choose the $d$-variate 
elliptically contoured $\alpha$-stable distribution, where no analytical formula
for the \pdf is known, and which generalises the normal law.  These distributions 
have heavy tails and are often used for modelling
financial data \citep{Nolan_stable_distrib}.  
We access the distribution in the next Example~\ref{ex:biVarEllip} through its  \pcf, which is
known analytically \refeq{eq:motiv_ch2}.

\begin{example}[$\alpha$-stable distribution]
\label{ex:biVarEllip}
The \emph{\pcf} of a $d$-variate elliptically contoured $\alpha$-stable distribution is given by
\begin{equation}  \label{eq:motiv_ch2}
\vphi_{\vX}(\bt)=\exp \left( \ii \bkt{\vek{t}}{\vek{\mu}}  - 
    \bkt{\vek{t}}{\bC \vek{t}}^{\frk{\alpha}{2}}\right).
\end{equation}
We approximate $\vphi_{\vX}(\bt)$ as in \refeq{eq:pcf2}, but in the TT format \refeq{eq:TTRepresentation}
and \refeq{eq:tt}. The tolerance used in the AMEn algorithm is $10^{-9}$. Further, 
from the inversion theorem, the \pdf of $\vX$ on $\D{R}^d$ can be computed as in 
\refeq{eq:motiv_pdf_lr} via the FFT.
\end{example}

We start by computing the KLD between two $\alpha$-stable distributions for fixed $\alpha_1 = 2.0$, 
$\alpha_2=1.9$ (with $\vek{\mu}_1=\vek{\mu}_2=0$, $\bC_1=\bC_2=\bI$); the results are summarised in Table~\ref{table:alpha_KLD_n_d}.
\begin{center}
\begin{table}[t]
 \centering
  \caption{Computation of $D_{\text{KL}}(\alpha_1,\alpha_2)$ for between two $\alpha$-stable distributions 
 for $\alpha_1 = 2.0$, $\alpha_2=1.9$, and different $d$ and $n$.  The AMEn tolerance is $10^{-9}$. }
 \label{table:alpha_KLD_n_d} 
 \begin{tabular}{lllllllllll}
 \hline
$d$ &16&16&16 &16 &16 & 16 & 16 & 32& 32&32 \\ \hline 
$n$ &8&16&32&64 & 128& 256&512& 64& 128&256\\ \hline
$D_{\text{KL}}(2.0,1.9)$ &0.016&0.059&0.06&0.062 & 0.06& 0.06 & 0.06&0.09 &0.14 & 0.12 \\ 
comp.time [s]   &0.8&3&8.9&14& 22 & 61 & 207  &46& 100& 258 \\ 
max. TT rank       &40&57&79&79& 59 & 79 & 77  &80&78& 79 \\ 
memory, MB &1.8&7&34&54& 73 & 158& 538 &160&313& 626 \\ \hline
\end{tabular}
\end{table}
\end{center}

From Table~\ref{table:alpha_KLD_n_d}  one may see that $n=32$ (for $d=16$) is sufficient, and there is
no need to a take higher resolution $n$, the KLD value is (almost) not changing. 
One may also see that for $d=32$ one needs to take $n=256$ or higher, but a higher $n$ 
requires more memory.  

These values in the last column in Table~\ref{table:alpha_KLD_n_d}, namely $d=32$ and $n=256$,
can be used to illustrate the amount of data and computation which would be involved in
a---here impossible---full representation.  The values $d=32$ and $n=256$ mean that the amount of data 
in full storage mode would be $N= n^d = 265^{32} \approx 1.16\mrm{E}77$,
and assuming 8 bytes per entry, this would be ca.\ $1\mrm{E}78$ bytes.  
Compare this to the estimated number of hadrons in the universe ($1\mrm{E}80$), to see that
alone the storage of such an object is not possible in full mode, whereas in a 
TT-low-rank approximation it did not require more than ca.\ $626\mrm{MB}$, and 
fits on a laptop.  And as for the computation of the KLD 
shown in that table, assume that the computation of the logarithms per data point could 
be achieved at a rate of $1 \mrm{GHz}$.  Then the KLD computation in full mode would 
require ca.\ $1.2\mrm{E}68 \mrm{s}$, or more than $3\mrm{E}60$ years, and even with a 
perfect speed-up on a parallel super-computer with say 
$1,000,000$ processors, this would require still more than $3\mrm{E}54$ years; compare this
with the estimated age of the universe of ca.\ $1.4\mrm{E}10$ years.

Continuing our tests,
in Table~\ref{table:divergences-comput_alpha_KLD} the KLD $D_{\text{KL}}(\alpha_1,\alpha_2)$ between 
two $\alpha$-stable distributions for different pairs of $(\alpha_1,\alpha_2)$ and fixed 
$d=8$ and $n=64$ is computed. The mean and covariance matrices were taken $\vek{\mu}_1=\vek{\mu}_2=0$, $\bC_1=\bC_2=\bI$. The tolerance for the AMEn algorithm was $10^{-12}$.  These results demonstrate that a 
TT approximation is possible (although the TT ranks are not so small) for various values of 
the parameters $\alpha$ in \refeq{eq:motiv_ch2}.
\begin{center}
\begin{table}[t]
 \centering
  \caption{Computation of $D_{\text{KL}}(\alpha_1,\alpha_2)$ between two $\alpha$-stable 
 distributions for various $\alpha$ with fixed $d=8$ and $n=64$. AMEn tolerance is $10^{-12}$.  $\vek{\mu}_1=\vek{\mu}_2=0$, $\bC_1=\bC_2=\bI$.}
 \label{table:divergences-comput_alpha_KLD} 
 \begin{tabular}{lllllll}
 \hline
$(\alpha_1, \alpha_2)$ & $(2.0,0.5)$ &$(2.0,1.0)$ &$(2.0,1.5)$ & $(2.0,1.9)$& $(1.5,1.4)$& $(1.0,0.4)$\\ \hline
$D_{\text{KL}}(\alpha_1,\alpha_2)$ & 2.27 & 0.66 & 0.3 & 0.03  & 0.031 & 0.6\\
comp.time [s] & 8.4  & 7.8 & 7.5&  8.5& 11 & 8.7 \\ 
max. TT rank & 78 & 74& 76& 76 & 80 & 79\\ 
memory, MB  & 28.5 & 28.5 & 27.1& 28.5& 35 & 29.5\\ \hline
\end{tabular}
\end{table}
\end{center}
%

%

From the KLD we turn again to the computation of the Hellinger distance,
this time for the $d$-variate 
elliptically contoured $\alpha$-stable distribution.
Table~\ref{table:divergences-comput_alpha_sqHellinger} shows the Hellinger distance 
$D_{H}(\alpha_1,\alpha_2)$ computed for two  different $\alpha$-stable distributions 
with values of $\alpha=1.5$ and $\alpha= 0.9$ for different $d$ and $n$. The mean values and the covariances are the same $\vek{\mu}_1=\vek{\mu}_2=0$, $\bC_1=\bC_2=\bI$. The maximal
TT ranks and the computation times are comparable to the KLD case in Table~\ref{table:alpha_KLD_n_d}.

\begin{center}
\begin{table}[h]
 \centering
  \caption{Computation of $D_{H}(\alpha_1,\alpha_2)$ between two $\alpha$ -stable distributions for different $d$ and $n$. AMEn tolerance is $10^{-9}$.  $\vek{\mu}_1=\vek{\mu}_2=0$, $\bC_1=\bC_2=\bI$. }
 \label{table:divergences-comput_alpha_sqHellinger}  
 \begin{tabular}{lllllllllll}
 \hline
$d$ &16&16&16 &16 &16 & 16 & 32 & 32 & 32& 32 \\ \hline 
$n$ &8&16&32&64 & 128& 256& 16 & 32& 64& 128\\ \hline
$D_{H}(1.5,0.9)$ &0.218& 0.223& 0.223& 0.223& 0.219& 0.223 &0.180 & 0.176& 0.175& 0.176  \\
comp.time [s]   &2.8& 3.7& 7.5&19  &53& 156&11& 21& 62& 117 \\
max. TT rank & 79& 76 & 76&76& 79& 76 &75& 71 & 75& 74     \\
memory, MB &7.7& 17& 34& 71&145& 283& 34& 66& 144& 285 \\
\end{tabular}
\end{table}
\end{center}

To show the influence of the TT (AMEn) tolerance, in
Table~\ref{table:Hellinger_conv} shows the $D_{H}$ distance computed with different 
TT (AMEn) tolerances. 
Additionally, the maximal tensor rank (there are $d$ ranks in total), the computing times, 
and the required storage cost are provided.

\begin{center}
\begin{table}[h]
 \centering
  \caption{Computation of $D_{H}(\alpha_1,\alpha_2)$ between two $\alpha$-stable distributions ($\alpha=1.5$ and $\alpha=0.9$) for different AMEn tolerances. $n=128$, $d=32$, $\vek{\mu}_1=\vek{\mu}_2=0$, $\bC_1=\bC_2=\bI$. }
 \label{table:Hellinger_conv} 
 \begin{tabular}{llllll}
 \hline
TT(AMEn) tolerance        & $10^{-7}$     & $10^{-8}$  & $10^{-9}$ & $10^{-10}$& $10^{-14}$ \\ \hline
$D_{H}(1.5,0.9)$  & 0.1645 & 0.1817 & 0.176 & 0.1761& 0.1802 \\ 
comp.time [s]& 43& 86& 103& 118& 241   \\
max. TT rank & 64& 75 & 75& 78& 77 \\
memory, MB &126& 255& 270& 307& 322   \\
\end{tabular}
\end{table}
\end{center}

\paragraph{A note about software.} Several tensor toolboxes developed for 
low-rank tensor calculus are available. 
The CP and Tucker decompositions are implemented in the Tensor Toolbox 
\citep{TTB_Dense, TTB_Software, TTB_CPOPT}, and in the Tensorlab
\citep{tensorlab3, DMV-SIAM2:00}. 
The TT and QTT tensor formats are implemented in \texttt{TT-toolbox} \citep{oseledets2011}.
The hierarchical Tucker tensor format is realised in the Hierarchical Tucker Toolbox 
\texttt{htucker toolbox} \citep{kressner2014algorithm}. 
For a more detailed overview, see \citep{GrasedyckKressnerTobler2013, khor-survey-2014}.
Almost all available implementations (e.g.\ TT-toolbox and htucker) are have been used to solve 
(stochastic) PDEs \citep{DolgovStat}, integral equations, linear systems in 
tensor format, or to perform arithmetic operations such as the scalar product, addition, etc.
But, to the best of our knowledge, we do not know any attempts for computing the KLD and other 
divergences or distances, or the entropy of high-dimensional probability distributions.
The Matlab files from our computations are freely available at 
\url{https://github.com/litvinen/Divergence_in_tensors.git}

\section{Conclusion} \label{S:concl}
The task considered here was the numerical computation of characterising statistics of 
high-dimensional \pdfs, as well as their divergences and distances, 
where the \pdf in the numerical implementation was assumed discretised on some regular grid.  
Even for moderate dimension $d$, the full storage and computation with such objects become very quickly 
infeasible.  

We have demonstrated that high-dimensional \pdfs, 
\pcfs, and some functions of them
can be approximated and represented in a low-rank tensor data format.  
Utilisation of low-rank tensor techniques helps to reduce the computational complexity 
and the storage cost from exponential $\C{O}(n^d)$ to linear in the dimension $d$, e.g.\ 
$\mathcal{O}\left(d n r^2\right )$ for the TT format. Here $n$ is the number of discretisation 
points in one direction, $r\ll n$ is the maximal tensor rank, and $d$ the problem dimension.
The particular data format is rather unimportant, 
any of the  well-known tensor formats (CP, Tucker, hierarchical Tucker, tensor-train (TT)) can be used, 
and we used the TT data format.  Much of the presentation and in fact the central train
of discussion and thought is actually independent of the actual representation.

In the beginning in \refS{intro} it was motivated through three possible ways how one may 
arrive at such a representation of the  \pdf.  One was if the \pdf was given in some approximate
analytical form, e.g.\  like a function tensor product of lower-dimensional \pdfs with a 
product measure, or from an analogous representation of the \pcf and subsequent use of the 
Fourier transform, or from a low-rank functional representation of a high-dimensional 
RV, again via its \pcf.
The theoretical underpinnings of the relation between \pdfs and \pcfs as well as their
properties were recalled in \refS{theory}, as they are important to be preserved in the
discrete approximation.  This also introduced the concepts of the convolution and of
the point-wise multiplication Hadamard algebra, concepts which become especially important if 
one wants to  characterise sums of independent RVs or mixture models,
a topic we did not touch on for the sake of brevity but which follows very naturally from
the developments here.  Especially the Hadamard algebra is also
important for the algorithms to compute various point-wise functions in the sparse formats.
The  \refS{theory}, as well as the following \refS{statistics} and   \refS{algs} are
actually completely independent of any particular discretisation and representation of the data

Some statistics, divergences, and distance measures were collected in \refS{statistics}
together with the abstract discrete expressions on how to compute them numerically,
independent of any particular numerical representation.  
To demonstrate our idea, one of the easiest tensor formats---the CP tensor format---was described
first in  \refS{tensor-rep}.  In the numerical part, we used the tensor-train (TT) format,
which is also sketched in \refS{tensor-rep}, together with how to implement the operations 
of the Hadamard algebra.
As \refS{statistics}
shows, point-wise functions are required in order to compute the desired statistics,
and some algorithms to actually perform this in an algebra were collected in  \refS{algs}.
These were originally developed for matrix algebra algorithms \citep{NHigham}, but they
work just as well in any other associative C$^*$-algebra.  In using such algorithms,
one assumes that the ranks of the involved tensors are not increasing strongly 
during iterations and linear algebra operations. 

In the numerical computations in \refS{Numerics}, the first example 
is one where the analytic answer was known, and this validates the approach
and shows the accuracy of the low-rank computation in the computation
of KL and Hellinger distances, which were
taken as examples of characterising functionals resp.\ statistics. 
As a more taxing problem, we then took elliptically contoured $\alpha$-stable distributions 
to evaluate the KLD and Hellinger distances between them.  For these distributions,
the \pdf is not known analytically, but the \pcf is, which we took as our
starting point. 
The AMEn TT-Cross algorithm \cite{ dolgov2015polynomial, DolgovLitvLiu19} was 
used to compute the low-rank approximations, and the \pdfs were then computed via FFT.
 This also nicely demonstrates the smooth integration of the FFT
into low-rank tensor formats, and
made it subsequently possible to compute the required quantities.

In total, these examples showed the viability of this concept, namely that it is possible to
numerically operate on discretised versions of high-dimensional distributions
with a reasonable computational expense, whereas in a full format the computations
would not have been feasible at all.
All that we required was that the data in its discretised form can be considered
as an element of a commutative C$^*$-algebra with an inner product,
where the algebra operations only have to be numerically computed in an
approximative fashion.
Such an algebra is isomorphic to a commutative sub-algebra
of the usual matrix algebra, allowing the use of matrix algorithms.





%
%
%
%
%
%
%
%
%



\section*{Acknowledgments}
The research reported in this publication was partly supported by funding 
from the Alexander von Humboldt Foundation (AvH), the Deutsche Forschungsgemeinschaft (DFG),
and a Gay-Lussac Humboldt Prize.
We also would like to thank Sergey Dolgov (University of Bath, UK) for his assistance with the TT-toolbox.



\begin{footnotesize}
\providecommand{\bysame}{\leavevmode\hbox to3em{\hrulefill}\thinspace}
\providecommand{\MR}{\relax\ifhmode\unskip\space\fi MR }
\providecommand{\MRhref}[2]{%
  \href{http://www.ams.org/mathscinet-getitem?mr=#1}{#2}
}
\providecommand{\href}[2]{#2}

\end{footnotesize}





\end{document}